\documentstyle{amsppt}
\magnification=1200
\NoBlackBoxes
\hsize=6.5truein
\vsize=8.9truein
\topmatter
\title Arithmetic Mixed Sheaves
\endtitle
\author Morihiko Saito
\endauthor
\affil RIMS Kyoto University, Kyoto 606-8502 Japan \endaffil
\endtopmatter
\tolerance=1000
\baselineskip=12pt
\def\scirc{\raise.2ex\hbox{${\scriptstyle\circ}$}}
\def\ssbull{\raise.2ex\hbox{${\scriptscriptstyle\bullet}$}}
\def\div{\text{\rm div}}
\def\dec{\text{\rm dec}}
\def\ind{\text{\rm ind}}
\def\hom{\text{\rm hom}}

\def\IC{\hbox{{\rm IC}}}
\def\CH{\hbox{{\rm CH}}}
\def\Spec{\hbox{{\rm Spec}}\,}
\def\Perv{\hbox{{\rm Perv}}}
\def\Ext{\hbox{{\rm Ext}}}
\def\Gr{\text{{\rm Gr}}}
\def\Im{\hbox{{\rm Im}}}
\def\Ker{\hbox{{\rm Ker}}}
\def\Coker{\hbox{{\rm Coker}}}
\def\Hom{\hbox{{\rm Hom}}}

\def\Gal{\hbox{{\rm Gal}}}
\def\Alb{\text{{\rm Alb}}}
\def\DR{\hbox{{\rm DR}}}
\def\MHM{\text{{\rm MHM}}}
\def\MHS{\text{{\rm MHS}}}
\def\supp{\hbox{{\rm supp}}\,}
\def\prim{\text{\rm prim}}
\def\Hdg{\text{\rm Hdg}}
\def\an{\text{\rm an}}
\def\Sh{\text{\rm Sh}}
\def\go{\text{\rm go}}
\def\red{\text{\rm red}}
\def\simto{\buildrel\sim\over\to}
\def\SameAuthor{\vrule height3pt depth-2.5pt width1cm}

\document
\noindent
We give a formalism of arithmetic mixed sheaves, including the case of
arithmetic mixed Hodge structures which are recently studied by P.
Griffiths, M.~Green ([22], [23]) and M.~Asakura [1] using other methods.
This notion became necessary for us to describe the image of Griffiths'
Abel-Jacobi map for a generic hypersurface (inspired by previous work of M.
Green [21] and C.~Voisin [45], [46]), and also to prove the following
variant of results of D.~Mumford [31] and A.~Roitman [34] suggested (and
proved in the case
$ \dim X = 2 $) by S.~Bloch in Exercise of Appendix to Lecture 1 of [7]:

Let
$ X $ be a smooth proper complex algebraic variety.
If there is a morphism of complex varieties
$ S \rightarrow X $ inducing a surjective morphism
$ \CH_{0}(S)_{\Bbb Q} \rightarrow \CH_{0}(X)_{\Bbb Q}, $ then
$ \Gamma (X, {\Omega }_{X}^{j}) = 0 $ for
$ j > \dim S $.

Actually, it turns out that the usual Hodge theory [12] is enough for this.
See Prop.~1.4 of Lect.~6 in [46].
But the attempt led us to the following formulation.

For a subfield
$ A $ of
$ {\Bbb R} $, and a subfield
$ k $ of
$ {\Bbb C} $ with finite transcendence degree, the category of
$ k $-finite mixed
$ A $-Hodge structures
$ \MHS(A)_{\langle k \rangle} $ is defined to be the inductive limit of the
category of admissible variations of mixed
$ A $-Hodge structures on
$ S = \Spec R $ over
$ k_{R} := \overline{k} \cap R $, with
$ R $ running over finitely generated smooth
$ k $-subalgebras of
$ {\Bbb C} $.
(Here an admissible variation on
$ S $ means an admissible variation of mixed
$ A $-Hodge structures on
$ S_{\Bbb C} := S\otimes _{k_{R}}{\Bbb C} $ in the sense of [27], [44]
such that the Hodge filtration, the connection and the polarization are
defined on
$ S/k_{R} $.)
More generally, for a complex algebraic variety
$ X $, there is a
$ k $-subalgebra
$ R $ of
$ {\Bbb C} $ as above such that
$ X $ is defined over
$ R $, and the category
$ \MHM(X,A)_{\langle k \rangle} $ of
$ k $-finite mixed Hodge Modules on
$ X $ is defined similarly.
See (2.1).
These can be extended to the mixed sheaves in the sense of [37] where
$ l $-adic sheaves can also be included, and the main theorems of this
paper hold in the generalized situation.
(In fact, they apply even to the subcategory consisting of the objects of
geometric origin.)

A similar category was defined in [37],(1.9) for a subfield
$ K $ of
$ {\Bbb C} $ finitely generated over
$ k $, by assuming that the fraction field of
$ R $ is
$ K $ in the above definition.
(This was inspired by [7], p.~1.20.)
Then it is enough to take further the inductive limit over
$ K $ in order to get the above category of arithmetic mixed Hodge Modules
(or sheaves, more generally).
See also 6.1.7 of [6] and 1.1 of Lect.~6 in [46].
The origin of the idea can be found in p.~86 and pp.~98--99 of [41].

The category of arithmetic mixed Hodge structures is closer to a
conjectural category of mixed motives of Beilinson [4] than the category of
graded-polarizable mixed
$ A $-Hodge structures
$ \MHS(A) $ in the usual sense [12], but is different from both (at least
if we do not restrict to the subcategory consisting of the objects of
geometric origin).
There is a natural functor
$ \iota : \MHS(X,A)_{\langle k \rangle} \rightarrow \MHS(X,A), $ but
this is not fully faithful.
See (2.5) (ii).
We have a constant object
$ A_{\langle k \rangle} $ in
$ \MHS(A)_{\langle k \rangle} $ represented by a constant variation of
Hodge structure of type
$ (0,0) $ on
$ S $.
For a complex variety
$ X $, we can define the cohomology
$ H^{j}(X, A_{\langle k \rangle}) $ in
$ \MHS(A)_{\langle k \rangle} $ in a compatible way with the functor
$ \iota $, and show the following:

\proclaim{{\bf 0.1.~Theorem}}
Let
$ X $ be a smooth proper complex variety of dimension
$ n $.
If
$ \Gamma (X, {\Omega }_{X}^{j}) \ne 0 $, then
$ \Ext^{j}(A_{\langle k \rangle}, H^{2n-j}(X,A_{\langle k \rangle}(n))) $
is an infinite dimensional vector space over
$ A \cap \overline{k} $, where
$ \Ext^{j} $ is taken in
$ \MHS(A)_{\langle k \rangle} $.
See (4.2).
\endproclaim

So we get a positive answer to the problem of showing the nonvanishing of
the higher extension groups in a variant of the category of mixed Hodge
structures.
See also [1], [22], [23], [48], etc.
The proof uses an argument similar to the passage from a normal function
(which is an element of an extension group) to the corresponding Hodge cycle
as in [38], [39].
The infinite dimensionality comes from the infinite transcendence degree of
$ {\Bbb C} $.
Note that (0.1) and the other main theorems hold for the arithmetic mixed
sheaves and also for the objects of geometric origin, where
$ A $ is assumed to be
$ {\Bbb Q} $.

For the relation with the cycle map, we define an analogue of Deligne
cohomology by
$$
\aligned
{H}_{\Cal D}^{i}(X,A_{\langle k \rangle}(j))= \Ext^{i}
&(A_{\langle k\rangle}, (a_{X})_{*}{a}_{X}^{*}A_{\langle k \rangle}(j))
\\
(=\Ext^{i}
&({a}_{X}^{*}A_{\langle k\rangle},{a}_{X}^{*}
A_{\langle k\rangle}(j))),
\endaligned
$$
where
$ a_{X} : X \rightarrow \Spec {\Bbb C} $ denotes the structure morphism,
and the extension group is taken in the derived category of
$ \MHS(A)_{\langle k \rangle} $ or
$ \MHM(X,A)_{\langle k \rangle} $.
By the decomposition theorem, we have the Leray filtration
$ L $ on
$ {H}_{\Cal D}^{i}(X,A_{\langle k \rangle}(j)) $ (associated with the Leray
spectral sequence) such that
$$
{\Gr}_{L}^{r}{H}_{\Cal D}^{i}(X,A_{\langle k \rangle}(j)) =
\Ext^{r}(A_{\langle k \rangle}, H^{i-r}(X, A_{\langle k \rangle}(j)).
$$
See (3.3).
For
$ A = {\Bbb Q} $, we have a cycle map
$$
cl : \CH^{p}(X)_{\Bbb Q} \rightarrow {H}_{\Cal D}^{2p}(X,{\Bbb
Q}_{\langle k \rangle}(p)),
$$
and the Chow group
$ \CH^{p}(X)_{\Bbb Q} $ has the induced filtration
$ L $.
(This can be generalized to Bloch's higher Chow groups [8] due to
[37], (8.3).)
It is expected that this filtration coincides with a conjectural filtration
of A.~Beilinson [4] and S.~Bloch.

By definition,
$ L^{1}\CH^{p}(X)_{\Bbb Q} $ consists of homologically equivalent to zero
cycles, and
$ L^{2}\CH^{p}(X)_{\Bbb Q} $ is contained in the kernel of Griffiths'
Abel-Jacobi map, but it is not clear if the last two coincide in general
(except when
$ p = 1 $ or
$ \dim X $).
The graded piece of the cycle map
$$
{\Gr}_{L}^{r}cl : {\Gr}_{L}^{r}\CH^{p}(X)_{\Bbb Q} \rightarrow
\Ext^{r}({\Bbb Q}_{\langle k \rangle}, H^{2p-r}(X, {\Bbb Q}_{\langle k\rangle}
(p)))
$$
for
$ r \ge 2 $ is sometimes called the higher Abel-Jacobi map.
See also [1], [22], [23], [48], etc.
By the same argument as in [40,~II], (3.3), we can verify that the
filtration
$ L $ is preserved by the action of a correspondence, and the action on the
graded pieces depends only on the cohomology class of the correspondence.
See (3.5).
Using Murre's result [32], [33], we can prove

\proclaim{{\bf 0.2.~Theorem}}
Let
$ X $ be a smooth proper variety of dimension
$ n $.
Then the image of the second Abel-Jacobi map
$ {\Gr}_{L}^{2}cl $ for
$ p = n $ is an infinite dimensional
$ {\Bbb Q} $-vector space if
$ \Gamma (X,{\Omega }_{X}^{2}) \ne 0 $.
See (4.4).
\endproclaim

In the surface case this is analogous to a recent result of C.~Voisin [48],
where she uses an explicit description of a second Abel-Jacobi map due to
M.~Green [22].
For a zero cycle
$ \zeta $ on a surface
$ S $, it is given by a composition of the extension classes associated
with the inclusions
$ |\zeta | \rightarrow C $ and
$ C \rightarrow S $, where
$ C $ is a curve containing
$ |\zeta | $.
But how to eliminate the cohomology of
$ C $ in the composition of the two extensions is a difficult problem,
and the relation between the two second Abel-Jacobi maps is not trivial.
(See also [1] for the nonvanishing of the second Abel-Jacobi map in the
case
$ X $ is the self-product of a curve defined over
$ \overline{\Bbb Q} $ such that its genus is
$ > 1 $ and the rank of the N\'eron-Severi group of the product over
$ \overline{\Bbb Q} $ is three.)

Note that the nonvanishing of
$ {\Gr}_{L}^{2}cl $ is rather trivial in the case
$ \Gamma (X,{\Omega }_{X}^{1}) = 0 $ and
$ X $ is defined over
$ k $ (i.e. there exists a
$ k $-variety
$ X_{k} $ such that
$ X = X_{k}\otimes _{k}{\Bbb C}) $,
because it is enough to consider a cycle of the form
$ [\Delta _{k}] - [X_{k}]\times \zeta _{k} $ where
$ \Delta _{k} $ is the diagonal of
$ X_{k} $ and
$ \zeta _{k} $ is a zero-cycle of degree
$ 1 $ on
$ X_{k} $ with rational coefficients.
(If
$ X $ is not defined over
$ k $,
we have to take a finitely generated
$ k $-subalgebra
$ R $ of
$ {\Bbb C} $ such that
$ X $ is defined over
$ R $,
and then carry out the argument relative to
$ R $.)

The nonvanishing of
$ {\Gr}_{L}^{r}cl $ for
$ p = n $ would hold under the assumption
$ \Gamma (X, {\Omega }_{X}^{r}) \ne 0 $, if the Chow-K\"unneth
decomposition in the sense of Murre [33] holds.
Note that the relation with Murre's result was suggested in [40,~II],
(3.4).
In the case
$ X $ is a smooth hypersurface of
$ {\Bbb P}^{n+1} $ defined over
$ k $, we can show that the image of
$ {\Gr}_{L}^{n}cl $ is infinite dimensional if
$ \Gamma (X,{\Omega }_{X}^{n}) \ne 0 $ (i.e.
if
$ \deg X \ge n + 2) $.
Similar arguments apply to cycles of lower codimensions.
See Remark after (4.4).
We can prove an analogue of (0.2) for the higher Chow groups.
See (5.2).

As to the image of Griffiths' Abel-Jacobi map, we have a natural morphism
of
$ {H}_{\Cal D}^{2p}(X,{\Bbb Q}_{\langle k \rangle}(p)) $ to the usual
Deligne cohomology
$ {H}_{\Cal D}^{2p}(X,{\Bbb Q}(p)) $, and its image is called the
$ k $-finite part of the Deligne cohomology, and is denoted by
$ {H}_{\Cal D}^{2p}(X,{\Bbb Q}(p))^{\langle k \rangle} $.
The
$ k $-finite part
$ J^{p}(X)_{\Bbb Q}^{\langle k \rangle } $ of Griffiths' intermediate
Jacobian tensored with
$ {\Bbb Q} $ is defined by intersecting it with
$ {H}_{\Cal D}^{2p}(X,{\Bbb Q}(p))^{\langle k \rangle} $.
(It coincides with the image of the extension group defined in the
category of arithmetic mixed Hodge structures.)
The Abel-Jacobi map induces a morphism of
$ L^{1}\CH^{p}(X)_{\Bbb Q} $ to
$ J^{p}(X)_{\Bbb Q}^{\langle k \rangle } $.
By an argument similar to [38], [39], we can show (see (4.1)):

\proclaim{{\bf 0.3.~Proposition}}
Let
$ X $ and
$ k $ be as in (0.1).
If
$ H^{2p-1}(X,{\Bbb Q}_{\langle k \rangle}) $ is global section-free in the
sense of (2.4), and the Hodge conjecture for codimension
$ p $ cycles holds for any smooth projective complex varieties defined over
$ k, $ then the cycle map
$ cl : \CH^{p}(X)_{\Bbb Q} \rightarrow {H}_{\Cal D}^{2p}(X,{\Bbb
Q}(p))^{\langle k \rangle} $ and the Abel-Jacobi map
$ L^{1}\CH^{p}(X)_{\Bbb Q} \rightarrow J^{p}(X)_{\Bbb Q}^{\langle k
\rangle } $ are surjective.
\endproclaim

The first assumption is satisfied if
$ X $ is a generic hypersurface, or more generally, if
$ X $ is a generic hypersurface of a smooth projective variety
$ Y $ defined over a subfield of
$ k $ such that
$ H^{2p-1}(Y,{\Bbb Q}) = 0 $.
If the first assumption of (0.3) is not satisfied, the cycle map
$ cl : \CH^{p}(X)_{\Bbb Q} \rightarrow {H}_{\Cal D}^{2p}(X,{\Bbb Q}
(p))^{\langle k \rangle} $ and the Abel-Jacobi map
$ L^{1}\CH^{p}(X)_{\Bbb Q} \rightarrow J^{p}(X)_{\Bbb Q}^{\langle k
\rangle } $ cannot be surjective in general, because
$ \Ext^{1} $ in
$ \MHS(k,A) $ is too big.
See Remark (ii) after (4.1).
(One possibility of getting smaller extension groups is to restrict to
objects of geometric origin.
See [40].)
We have an analogue of (0.3) for higher Chow groups.
See (5.11).

One of the main problems in this theory is the injectivity of the cycle
map defined before (0.2).
We will reduce this problem to the injectivity of certain Abel-Jacobi map
for varieties over number fields.
See also Example (3.7).

Let
$ Y $ be a smooth variety over an algebraic number field
$ k $, where
$ k $ is assumed to be algebraically closed in the function field of
$ Y $.
We consider the Abel-Jacobi map over
$ k $:
$$
L^{1}\CH^{p}(Y)_{\Bbb Q} \rightarrow J^{p}(Y/k)_{\Bbb Q} :=
\Ext^{1}({\Bbb Q}_{k}, H^{2p-1}(Y/k,{\Bbb Q}(p)))
$$
where
$ L^{1}\CH^{p}(Y)_{\Bbb Q} $ is the subgroup of cohomologically
equivalent to zero cycles (with respect to
$ H^{\ssbull }(Y\otimes _{k }\overline{k}, {\Bbb Q}_{l}) $ or
$ H^{\ssbull }(Y({\Bbb C}),{\Bbb Q})) $, and
$ {\Bbb Q}_{k}, H^{2p-1}(Y/k,{\Bbb Q}(p)) $ are defined in
$ \MHS(k,A) $, the category of graded-polarizable mixed
$ {\Bbb Q} $-Hodge structures whose
$ {\Bbb C} $-part is defined over
$ k $.
See (1.3).
(In the case
$ k = $
$ {\Bbb C} $ and
$ p = 1 $, we can show that this Abel-Jacobi map is bijective.
See [40,~I], (3.4).)

For a smooth complex algebraic variety
$ X $, we say that
$ Y $ is a
$ k $-smooth (or
$ k $-smooth proper) model of
$ X $, if
$ Y $ is a smooth (or smooth proper)
$ k $-variety having a
$ k $-morphism to an integral
$ k $-variety whose geometric generic fiber over
$ \Spec {\Bbb C} $ is isomorphic to
$ X $.

\proclaim{{\bf 0.4.~Theorem}}
Let
$ X, k $ be as in (0.1), and assume
$ k $ is a number field.
Then the cycle map
$ \CH^{p}(X)_{\Bbb Q} \rightarrow {H}_{\Cal D}^{2p}(X,{\Bbb Q}_{\langle
k \rangle}(p)) $ is injective if the above Abel-Jacobi map over
$ k $ is injective for any
$ k $-smooth models
$ Y $ of
$ X $.
In the case
$ p = 2 $, the last assumption is reduced to the same injectivity for any
$ k $-smooth proper models
$ Y $ of
$ X $, and moreover,
$ Y $ can be replaced with a smooth projective
$ k $-variety which is birational to
$ Y $.
See (4.6-7).
\endproclaim

This shows that the second hypothesis in a theorem of [1] is not necessary.
It is expected that the above Abel-Jacobi map for smooth projective
varieties over number fields would be injective if the target is replaced
by the motivic extension group.
See also [4], [25].
However, it is not clear whether the morphism of the motivic extension
group to the above extension group should be injective, because it is
closely related to the full faithfulness of the forgetful functor from the
category of mixed motives to that of mixed Hodge structures, or of systems
of realizations.
See [40,~I].
So it is better to use the systems of realizations here.

If the conclusion of (0.4) is true, and the K\"unneth decomposition in the
Chow group holds for the smooth proper variety
$ X, $ then we can verify that
$ L $ coincides with the filtration of Murre [33], using the same argument
as in [40,~II], (3.3).
See (4.9) below.
We can also show that Bloch's conjecture in [7] on the converse of Mumford
's result [31] can be reduced to the last hypothesis of (0.4) (although the
latter may be more difficult than the former).
See Remark (i) after (4.8).
We can prove an analogue of (0.4) for the higher Chow groups.
See (5.3).
This gives evidence for a conjecture of C.~Voisin [47] on the countability
of
$ \CH^{2}(X,1)_{\Bbb Q} $ modulo decomposable higher cycles.
See (5.10).

In \S 1, we introduce the notion of mixed Hodge structure whose
$ {\Bbb C} $-part is defined over
$ k $, and, more generally, that of mixed sheaf.
In \S 2, we study
$ k $-finite mixed Hodge Modules (and also mixed sheaves), and define the
cycle map in \S 3.
The proofs of the above assertions in the generalized situation are given
in \S 4.
We show analogues of (0.2) and (0.4) for the higher Chow groups in \S 5.

\bigskip\bigskip\centerline{{\bf 1.~Mixed Sheaves Defined over a Subfield}}

\bigskip
\noindent
{\bf 1.1.}
Let
$ A $ be a subfield of
$ {\Bbb R} $, and
$ k $ a subfield of
$ {\Bbb C} $.
For varieties
$ X $ over
$ k $, let
$ \MHM(X/k,A) $ denote the categories of mixed Hodge Modules [36] on
$ X_{\Bbb C} := X\otimes _{k}{\Bbb C} $ such that the underlying filtered
$ {\Cal D} $-Modules
$ (M_{\Bbb C},F;W) $ on
$ X_{\Bbb C} $ are defined over
$ k $ (i.e., a filtered
$ {\Cal D} $-Module
$ (M_{k},F;W) $ on
$ X/k $ together with an isomorphism
$ (M_{k},F;W)\otimes _{k}{\Bbb C} = (M_{\Bbb C},F;W) $ is given).
We also assume that polarizations on the graded pieces of the weight
filtration are defined over
$ k $.
See (1.8, ii) in [37].

More generally, it is possible to consider the categories of mixed sheaves
$ {\Cal M}(X/k,A) $ for varieties
$ X $ over
$ k $ which satisfy the axioms of mixed sheaves.
See loc.~cit.
Here we have to choose an embedding
$ k \rightarrow {\Bbb C} $ for the later argument.
We assume furthermore that for
$ k \subset k' \subset {\Bbb C} $ and for a
$ k $-variety
$ X $, we have canonically the base change functor
$$
{\Cal M}(X/k,A) \rightarrow {\Cal M}(X\otimes _{k}k'/k',A),
\leqno(1.1.1)
$$
which is exact and compatible with dual, pull-back, direct image, external
product, etc.
We also assume that there exist forgetful functors
$$
{\Cal M}(X/k,A) \rightarrow \MHM(X_{\Bbb C},A)
\leqno(1.1.2)
$$
compatible with the forgetful functors to
$ \Perv(X_{\Bbb C},A) $.

We can consider, for example, the category consisting of
$ ((M_{k};F,W), (K,W), (K_{l},W)) $ where
$ (M_{k},F) $ is a regular holonomic filtered
$ {\Cal D} $-Module on
$ X/k $ endowed with a finite filtration
$ W $,
$ (K,W) $ is a filtered perverse sheaf on
$ {X}_{\Bbb C}^{\an} $ with
$ {\Bbb Q} $-coefficients, and
$ (K_{l},W) $ are filtered perverse
$ l $-adic sheaves on
$ X_{\overline{k}} := X\otimes _{k}\overline{k} $ with
$ {\Bbb Q}_{l} $-coefficients which are endowed with a continuous action of
the Galois group of
$ \overline{k}/k $ as in [39] (i.e. the action is lifted to perverse
sheaves with
$ {\Bbb Z}_{l} $-coefficients).
Furthermore, they are given comparison isomorphisms
$$
\DR((M_{k},W)\otimes _{k}{\Bbb C}) = (K,W)\otimes _{\Bbb Q}
{\Bbb C},\quad \varepsilon ^{*}{i}_{\overline{k}}^{*}(K_{l},W) =
(K,W)\otimes _{\Bbb Q}{\Bbb Q}_{l}.
\leqno(1.1.3)
$$
Here
$ A = {\Bbb Q} $,
$ \overline{k} $ is the algebraic closure of
$ k $ in
$ {\Bbb C} $, and
$ i_{\overline{k}}:X_{\Bbb C} \rightarrow X_{\overline{k}} $ is the
canonical morphism.
(See [6] for
$ \varepsilon ^{*} $.)
It is also possible to consider the perverse sheaves
$ K_{\sigma} $ on
$ (X \otimes _{k,\sigma} {\Bbb C})^{\an} $ for any embeddings
$ \sigma : k \rightarrow {\Bbb C} $ so that we get the systems of
realizations as in [37], (1.8).
See [14], [15], [25], etc. for the case
$ X = \Spec k $.

\medskip
\noindent
{\it Remark.} If
$ X $ is a smooth complex variety, a mixed
$ A $-Hodge Module consists of
$ ((M,F;W), (K,W), \alpha ) $ where
$ (M,F) $ is a filtered regular holonomic
$ {\Cal D}_{X} $-Module,
$ K $ is a perverse sheaf [6] with
$ A $-coefficients on
$ X^{\an} $,
$ W $ is a finite increasing filtration on
$ M, K $ (called the weight filtration), and
$ \alpha : \DR(M,W) = (K,W)\otimes _{A}{\Bbb C} $ is an isomorphism of
filtered perverse sheaves with
$ {\Bbb C} $-coefficients on
$ X^{\an} $.
They satisfy several good conditions.
See [36].
Note that the condition of mixed Hodge Module on
$ {\Cal M} = ((M,F;W), (K,W), \alpha ) $ is Zariski-local, and we have
locally the following: If
$ g $ is a function on
$ X $ such that the restriction of
$ K $ to the complement
$ U $ of
$ Y := g^{-1}(0)_{\red} $ is a local system up to a shift of complex, then
$ {\Cal M} $ is a mixed Hodge Module if and only if the following four
conditions are satisfied :

\noindent
(a) The restriction of
$ {\Cal M} $ to
$ X \setminus Y $ is an admissible variation of mixed Hodge structure in
the sense of [27], [44].

\noindent
(b) The three filtrations
$ F, V, W $ on
$ (i_{g})_{*}M $ (the direct image as a
$ {\Cal D} $-Module) are compatible, where
$ i_{g} $ is the closed embedding by the graph of
$ g $.

\noindent
(c) The relative monodromy filtration exists on
$ \psi _{g}(K,W) $ and
$ \varphi _{g,1}(K,W) $.

\noindent
(d)
$ \psi _{g}{\Cal M} $ and
$ \varphi _{g,1}{\Cal M} $ are mixed Hodge Modules.

\noindent
See [34] for the details.

If
$ X $ is a singular complex variety, the underlying filtered
$ {\Cal D} $-Module of a mixed
$ A $-Hodge Module is defined by using closed embeddings
$ U \rightarrow V $ where
$ U $ is an open subvariety of
$ X $ and
$ V $ is smooth.
See for example [35], 2.1.20.

\proclaim{{\bf 1.2.~Theorem}}
With the above notation, the category
$ {\Cal M}(X/k,A) $ is an
$ (A \cap k) $-linear abelian category, and every morphism is strictly
compatible with
$ (F,W) $ on
$ M_{k} $.
The bounded derived categories
$ D^{b}{\Cal M}(X/k,A) $ are stable by standard functors like
$ f_{*}, f_{!}, f^{*}, f^{!}, $ etc. in a compatible way with (1.1.2).
\endproclaim

\demo\nofrills {Proof.\usualspace}
This follows from [37].
In the case
$ {\Cal M}(X/k,A) = \MHM(X_{\Bbb C},A) $, the first assertion follows
from [35], 5.1.14, and the last from 4.3 and 4.4 in loc.~cit.
\enddemo

\noindent
{\it Remark.} For the construction of the direct images
$ f_{*}, f_{!} $, it is enough to construct the cohomological direct images
$ H^{i}f_{*}, H^{i}f_{!} $ (thanks to the strictness of
$ (F,W)) $.
Indeed, if
$ X $ is quasi-projective,
$ f_{*}, f_{!} $ are defined by taking two sets of open coverings of
$ X $ associated with general hyperplane sections, and using the
combination of co-Cech and Cech complexes together with Artin's vanishing
theorem for the (perverse) cohomological direct images by an affine
morphism.
See [5], [6].

For the pull-back by a closed embedding
$ i : X \rightarrow Y $, it is enough to show an equivalence of categories
$$
i_{*} : D^{b}\MHM(X/k,A) \rightarrow D_{X}^{b}\MHM(Y/k,A),
$$
where the target is the full subcategory of
$ D^{b}\MHM(Y/k,A) $ defined by the condition :
$ \supp H^{j}M \subset X $ for any
$ j $.
This is reduced to the case where
$ X $ is a divisor defined by a function
$ g $.
Then the assertion is verified by using the functor
$ \xi _{g} $ in [36], 2.22.
This is inspired by [5].

\medskip
\noindent
{\bf 1.3.}
Let
$ \MHS(k,A) = \MHM(\Spec k/k,A) $ (where we put
$ F^{p} = F_{-p}) $.
This is the category of graded-polarizable mixed
$ A $-Hodge structures [12] whose
$ {\Bbb C} $-part is defined over
$ k $ (and polarizations are assumed to be defined over
$ k $ as above).
An object
$ H $ of
$ \MHS(k,A) $ consists of
$ ((H_{k},F;W), (H_{A},W), (H_{\Bbb C},W); \alpha _{k}, \alpha _{A}) $,
where
$ H_{\Lambda } $ is a filtered (or bifiltered)
$ \Lambda $-vector spaces for
$ \Lambda = k, A, {\Bbb C} $ with comparison isomorphisms
$ \alpha _{k} : H_{k}\otimes _{k}{\Bbb C} \rightarrow H_{\Bbb C} $,
$ \alpha _{A} : H_{A}\otimes _{A}{\Bbb C} \rightarrow H_{\Bbb C} $
compatible with the filtration
$ W $.

For
$ j \in {\Bbb Z} $, we define
$ A_{k}(j) \in \MHS(k,A) $ by
$ H_{k} = k $ and
$ H_{A} = (2\pi i)^{j}A \subset H_{\Bbb C} = {\Bbb C} $, where
$ \alpha _{k} $ and
$ \alpha _{A} $ are natural isomorphisms, and
$ {\Gr}_{F}^{i} = 0 $ for
$ i \ne -j $ and
$ {\Gr}_{i}^{W} = 0 $ for
$ i \ne -2j $.
Let
$ A_{k} = A_{k}(0) $.
We have also
$ A_{k}(j) \in {\Cal M}(k,A) $ in general.
See [37].

For a
$ k $-variety
$ X $ with structure morphism
$ a_{X/k} : X \rightarrow \Spec k, $ let
$$
A_{X/k}(j) = {a}_{X/k}^{*}A_{k}(j),\quad H^{j}(X/k, A(j)) =
H^{j}(a_{X/k})_{*}A_{X/k}(j).
\leqno(1.3.1)
$$

For a smooth
$ k $-variety
$ X $, we have a cycle map
$$
cl_{k} : \CH^{p}(X)_{\Bbb Q} \rightarrow \Ext^{2p}(A_{X/k},
A_{X/k}(p))=\Ext^{2p}(A_{k}, (a_{X/k})_{*}A_{X/k}(p)),
\leqno(1.3.2)
$$
where the extension groups are taken in the derived categories of
$ {\Cal M}(X/k,A) $ or
$ {\Cal M}(k,A) $.
If a cycle
$ \zeta $ is represented by an irreducible closed subvariety
$ Z $, then
$ cl(\zeta ) $ is represented by the composition of
$$
A_{X/k} \rightarrow A_{Z/k} \rightarrow \IC_{Z/k}A[-d_{Z/k}]
$$
with its dual, by using the dualities
$$
{\Bbb D}(A_{X/k}) = A_{X/k}(d_{X/k})[2d_{X/k}],\quad {\Bbb D}(\IC_{Z/k}A) =
\IC_{Z/k}A (d_{Z/k}),
$$
where
$ \IC_{Z/k}A $ is the intersection complex, and
$ d_{X/k} = \dim X/k $.
See [36], (4.5.15).

This cycle map is compatible with the pushdown and the pull-back of cycles
(where we assume that a morphism is proper for the pushdown.)
See [40,~II].
We define
$$
L^{1}\CH^{p}(X)_{\Bbb Q} = \Ker(\CH^{p}(X)_{\Bbb Q} \rightarrow
\Hom({\Bbb Q}_{k}, H^{2p}(X/k,{\Bbb Q}(p))))
\leqno(1.3.3)
$$
where the last morphism is the composition of the cycle map with the
natural morphism to
$ \Hom({\Bbb Q}_{k}, H^{2p}(X/k,{\Bbb Q}(p))) $.
Then the cycle map induces the Abel-Jacobi map over
$ k $
$$
L^{1}\CH^{p}(X)_{\Bbb Q} \rightarrow J^{p}(X/k)_{\Bbb Q} :=
\Ext^{1}({\Bbb Q}_{k}, H^{2p-1}(X/k,{\Bbb Q}(p)))
\leqno(1.3.4)
$$
by using the Leray spectral sequence
$$
{E}_{2}^{i,j} = \Ext^{i}({\Bbb Q}_{k}, H^{j}(X/k,{\Bbb Q}(p))) \Rightarrow
\Ext^{i+j}({\Bbb Q}_{X/k}, {\Bbb Q}_{X/k}),
$$
because the vanishing of
$ {E}_{2}^{i,j} $ for
$ i < 0 $ implies the injective morphisms
$$
{E}_{\infty }^{i,j} \rightarrow {E}_{2}^{i,j}\,\,\,(i \le 1).
\leqno(1.3.5)
$$

Similar assertions hold also for
$ {\Cal M}(X/k,A) $.
See [37].

\proclaim{{\bf 1.4.~Lemma}}
For
$ H_{i} \in \MHS(k,A) (i = 1, 2), $ let
$ H = {\Cal H}om(H_{1},H_{2}) \in \MHS(k,A) $, and
$$
C(H)^{0} = \Ker(F^{0}W_{0}H_{k}\oplus W_{0}H_{A} \rightarrow {\Gr}_{0}^{W}
H_{\Bbb C}),\quad C(H)^{1} = W_{-1}H_{\Bbb C},
$$
so that we get a complex
$ C(H) = [C(H)^{0} \rightarrow C(H)^{1}] $.
Then we have natural isomorphisms
$$
H^{0}(C(H)) = \Hom(H_{1},H_{2}),\,\,\,H^{1}(C(H)) = \Ext^{1}(H_{1},H_{2}),
$$
where
$ \Hom $ and
$ \Ext^{1} $ are taken in
$ \MHS(k,A) $.
\endproclaim

(The proof is left to the reader.
The argument is similar to [11], taking account of the semisimplicity of
$ {\Gr}_{0}^{W}H $.)

\proclaim{{\bf 1.5.~Corollary}}
For
$ i > 1 $,
$ \Ext^{i}(H_{1},H_{2}) = 0 $.
\endproclaim

\bigskip\bigskip\centerline{{\bf 2.~Limit of Mixed Sheaves}}

\bigskip
\noindent
{\bf 2.1.}
$ k $-{\it finite Mixed Sheaves.} Let
$ A, k $ be as in (1.1).
Let
$ K $ be a subfield of
$ {\Bbb C} $ containing
$ k $.
For a
$ K $-variety
$ X $, let
$ R $ be a finitely generated smooth
$ k $-subalgebra of
$ K $ such that
$ X $ is defined over
$ R $, i.e., there is an
$ R $-scheme
$ X_{R} $ with an isomorphism
$ X_{R}\otimes _{R}K = X $ (and
$ \Spec R $ is smooth over
$ k $).
Then, for a finitely generated smooth
$ k $-subalgebra
$ R' $ of
$ {\Bbb C} $ containing
$ R $, we put
$ X_{R'} = X_{R}\otimes _{R}R', S' = \Spec R', d_{R'} = \dim_{k} S' $, and
$ k_{R'} = \overline{k} \cap R' $, where
$ \overline{k} $ is the algebraic closure of
$ k $ in
$ {\Bbb C} $.
Note that
$ k_{R'} $ coincides with the algebraic closure of
$ k $ in the function field
$ k(S') $ of
$ S' $ (because
$ S' $ is normal) so that
$ S/k_{R'} $ is geometrically irreducible and
$ S'_{\Bbb C} := S'\otimes _{k_{R'}}{\Bbb C} $ is connected.

Let
$ R', R'' $ be as above such that
$ R' \subset R'' $ and
$ S'' := \Spec R'' $ is smooth over
$ S' $.
Let
$ {\Cal M}(X/k,A) $ be as in (1.1) (e.g.
$ {\Cal M}(X/k,A) = \MHM(X/k,A)) $.
Then we have an exact functor by the composition of pull-backs
$$
\aligned
{\Cal M}(X_{R'}/k_{R'},A)[-d_{R'}]
&\rightarrow {\Cal M}(X_{R'}\otimes_{k_{R'}}k_{R''}/k_{R''},A)[-d_{R'}]
\\
&\rightarrow {\Cal M}(X_{R''}/k_{R''},A)[-d_{R''}],
\endaligned
\leqno(2.1.1)
$$
where
$ {\Cal M}(X_{R'}/k_{R'},A)[-d_{R'}] $ denotes the category of mixed
sheaves
$ {\Cal M}(X_{R'}/k_{R'},A) $ shifted by
$ -d_{R'} $ in the derived category.
(Note that
$ R'\otimes _{k_{R'}}k_{R''} \rightarrow R'' $ is injective.)
Then
$ X_{\Bbb C}\,(:=X\otimes_{K}{\Bbb C} $) is identified with the
closed fiber of
$ X_{R'}\otimes _{k_{R'}}{\Bbb C} $ over the closed point of
$ S'_{\Bbb C} $ defined by the inclusion
$ R' \subset {\Bbb C} $, and we have an exact functor
$$
{\Cal M}(X_{R'}/k_{R'},A)[-d_{R'}] \rightarrow \MHM(X_{\Bbb C},A)
\leqno(2.1.2)
$$
compatible with the above functor by using (1.1.2) and restricting to the
fiber.
Let
$ I $ be the ordered set consisting of finitely generated smooth
$ k $-subalgebras
$ R' $ of
$ {\Bbb C} $ containing
$ R $, where
$ R' < R'' $ if and only if
$ R' \subset R'' $ and
$ R'' $ is smooth over
$ R' $.
We define
$$
{\Cal M}(X/K,A)_{\langle k \rangle} = \hbox{\rm ind\,lim}_{I}
{\Cal M}(X_{R'}/k_{R'},A)[-d _{R'}].
\leqno(2.1.3)
$$
This is independent of the choice of
$ R $ and
$ X_{R} $ at the beginning of this section.
An object of this category is represented by an object of
$ {\Cal M}(X_{R'}/k,A)[-d_{R'}] $ for some
$ R' $, and is called a
$ k $-finite mixed
$ A $-Hodge Module.
By (2.1.2) we have a natural functor
$$
\iota : {\Cal M}(X/K,A)_{\langle k \rangle} \rightarrow
\MHM(X_{{\Bbb C}},A).
\leqno(2.1.4)
$$
We have (shifted) cohomological functors
$ H^{i} $ to
$ {\Cal M}(X_{R'}/k_{R'},A)[-d_{R'}] $,
$ {\Cal M}(X/K,A)_{\langle k \rangle} $,
$ \MHM(X_{{\Bbb C}},A) $ from their derived categories such that the
restrictions of
$ H^{0} $ to these subcategories are the identity, and the
$ H^{i} $ are
compatible with the functors (2.1.1), (2.1.2) and (2.1.4).

In the case the mixed sheaves
$ {\Cal M}(X/k,A) $ are mixed Hodge Modules
$ \MHM(X/k,A) $ in (1.1), we will denote
$ {\Cal M}(X/K,A)_{\langle k \rangle} $ by
$ \MHM(X/K,A)_{\langle k \rangle} $.
If
$ K = {\Bbb C} $,
$ {\Cal M}(X/K,A)_{\langle k \rangle} $ will be denoted by
$ {\Cal M}(X,A)_{\langle k \rangle} $ to simplify the notation.

\medskip
\noindent
{\bf 2.2.}
{\it Remarks.} (i) If
$ K $ has finite transcendence degree, let
$ k' $ be the algebraic closure of
$ k $ in
$ K $, and
$ d_{K/k} $ the transcendence degree of
$ K $ over
$ k $.
Then we have
$$
{\Cal M}(X/K,A)_{\langle k \rangle}
= \hbox{\rm ind\,lim}_{R'} {\Cal M}(X_{R'}/k',A)[-d _{K/k}],
\leqno(2.2.1)
$$
where
$ R' $ runs over localizations of
$ R $ by one element, and
$ X_{R'} = X_{R}\otimes _{R}R' $.
(This is considered in [37], (1.9), see also [41], 2.12.)

(ii)
For a complex variety
$ X $, there exist a subfield
$ K $ of
$ {\Bbb C} $ which contains
$ k $ and is finitely generated over
$ k $, and a
$ K $-variety
$ X_{K} $ with an isomorphism
$ X_{K}\otimes _{K}{\Bbb C} = X $.
Then we have by definition
$$
{\Cal M}(X,A)_{\langle k \rangle} = \hbox{\rm ind\,lim}_{K'}
{\Cal M}(X_{K}\otimes _{K}K'/K',A)_{\langle k \rangle}[-d_{K'/k}],
\leqno(2.2.2)
$$
where
$ K' $ runs over subfields of
$ {\Bbb C} $ which contains
$ K $ and is finitely generated over
$ k $.

\medskip
\noindent
{\bf 2.3.}
$ k $-{\it finite Mixed Hodge structures.} With the notation of (1.1) and
(2.1), let
$$
{\Cal M}(K,A)_{\langle k \rangle}
= {\Cal M}(\Spec K/K,A)_{\langle k \rangle}.
\leqno(2.3.1)
$$
In the case
$ K = {\Bbb C} $, this will be denoted by
$ {\Cal M}(A)_{\langle k \rangle} $ to simplify the notation.
If the mixed sheaves
$ {\Cal M}(X/k,A) $ in (1.1) are mixed Hodge Modules
$ \MHM(X/k,A) $, it is denoted by
$ \MHS(K,A)_{\langle k \rangle} $.
If furthermore
$ k = {\Bbb C} $, it is denoted by
$ \MHS(A)_{\langle k \rangle} $, and is called the category of
$ k $-finite mixed
$ A $-Hodge structures.

In general, we have a natural functor
$$
\iota : {\Cal M}(K,A)_{\langle k \rangle} \rightarrow \MHS(A),
\leqno(2.3.2)
$$
where
$ \MHS(A) $ is the category of graded-polarizable mixed
$ A $-Hodge structures in the usual sense [12].

For a finitely generated smooth
$ k $-algebra
$ R $ of
$ {\Bbb C} $, put
$ S = \Spec R $ and
$ S_{\Bbb C} = \Spec R\otimes _{k_{R}}{\Bbb C} $.
Then admissible variations of mixed Hodge structures are mixed Hodge
Modules by [36], 3.27, and
$ \MHS(K,A)_{\langle k \rangle} $ is the inductive limit of the category
of admissible variations of mixed Hodge structures on
$ \Spec R/k_{R} $.

For
$ j \in {\Bbb Z} $, we have
$ A_{K,\langle k \rangle}(j) \in \MHS(K,A)_{\langle k \rangle} $ which is
represented by a constant variation of Hodge structure
$ A_{S/k_{R}}(j) $ of type
$ (-j,-j) $ on
$ S = \Spec R $.
It is denoted by
$ A_{\langle k \rangle}(j) $ if
$ K = \Spec {\Bbb C} $.
In general, we have
$ A_{K,\langle k \rangle}(j) \in {\Cal M}(K,A)_{\langle k \rangle},
A_{\langle k \rangle}(j)\in {\Cal M}(A)_{\langle k \rangle} $ similarly.

\medskip
\noindent
{\bf 2.4.~Definition.} We say that
$ H \in \MHS(K,A)_{\langle k \rangle} $ (or
$ {\Cal M}(K,A)_{\langle k \rangle} $)
is global section-free, if for any finitely generated smooth
$ k $-subalgebra
$ R $ of
$ K $ such that
$ H $ is represented by an admissible variation of mixed Hodge structure
(or a smooth mixed sheaf) on
$ S := \Spec R $, the underlying local system on
$ S_{\Bbb C} $ has no nonzero global section.

\medskip
\noindent
{\bf 2.5.}
{\it Remarks.} (i) When
$ k = \overline{\Bbb Q} $ and
$ K = {\Bbb C} $,
$ \MHS(A)_{\langle k \rangle} $ is called the category of arithmetic mixed
Hodge structures according to M.~Asakura [1], and is closely related to
recent work of P.~Griffiths and M.~Green [22], [23] (see also [48]).

If
$ H \in \MHS(A) $ is the image of
$ {\Cal M} \in \MHS(A)_{\langle k \rangle} $ by the functor
$ \iota $ in (2.3.2), then the underlying
$ {\Bbb C} $-part
$ H_{\Bbb C} $ has naturally an integrable connection
$ \nabla : H_{\Bbb C} \rightarrow H_{\Bbb C}
\otimes _{\Bbb C}{\Omega }_{{\Bbb C}/\overline{k}}^{1} $ which is the
limit of the connection
$ M_{R} \rightarrow M_{R}\otimes _{R}{\Omega }_{R/{k}_{R}}^{1} $, where
$ M_{R} = \Gamma (S,M_{S}) $ with
$ M_{S} $ the underlying
$ {\Cal D}_{S} $-Modules of the representatives of
$ {\Cal M} $.

So we can also consider the category
$ \MHS(A)_{\nabla /\overline{k}} $ of mixed Hodge structures whose
$ {\Bbb C} $-part has an integrable connection
$ \nabla $ over
$ \overline{k} $.
Then we can show that
$ \MHS(A)_{\langle k \rangle} $ is a full subcategory of
$ \MHS(A)_{\nabla /\overline{k}} $.
Indeed, a morphism
$ A \rightarrow \iota ({\Cal M}) $ in
$ \MHS(A)_{\nabla /\overline{k}} $ defines a section of
$ \Ker \,\nabla $ in
$ M_{R} $ for
$ R $ sufficiently large, and the corresponding global section of
$ H_{\Bbb C} $ on
$ S_{\Bbb C} $ belongs to the
$ A $-local system
$ H_{A} $ at the point corresponding to the inclusion
$ R \rightarrow {\Bbb C} $.
See also [1], [22], [23].

(ii) The canonical functor
$ {\Cal M}(A)_{\langle k \rangle} \rightarrow \MHS(A) $ is not fully
faithful, because the morphism of extension classes is not injective.
For example, let
$ S = \Spec R $ with
$ R = k[t,t^{-1}] $, and consider an admissible variation of mixed Hodge
structure on
$ S $ which is an extension of
$ {\Bbb Q}_{S} $ by
$ {\Bbb Q}_{S}(1) $, and which has an
$ {\Cal O}_{S} $-basis
$ e_{0}, e_{1} $ such that
$ t\partial _{t}e_{0} = e_{1}, t\partial _{t}e_{1} = 0 $,
$ F^{0} $ is generated by
$ e_{0} $ over
$ {\Cal O}_{S} $, and the
$ {\Bbb Q} $-lattice of the fiber at
$ t $ is generated by
$ e_{1} $ and
$ e_{0} + (2\pi i)^{-1}\log(t/c)e_{1} $, where
$ c \in {\Bbb C}^{*} $ corresponds to the point of
$ S_{\Bbb C} $ determined by the embedding
$ R \rightarrow {\Bbb C} $.

For another example, we can use the noninjectivity of
$ \Ext^{1}(A_{k},H) \rightarrow \Ext^{1}(A,{a}_{{\Bbb C}/{k}}^{*}H) $ for
$ H \in \MHS(A)_{k} $ (see (1.4)), where
$ {a}_{{\Bbb C}/{k}}^{*}H $ is the underlying graded-polarizable mixed
$ A $-Hodge structure of
$ H $.
Indeed, the morphism is factorized by the injective morphism
$ \Ext^{1}(A_{k},H) \rightarrow \Ext^{1}(A_{S/k},{a}_{S/{k}}^{*}H) $.

(iii) If the objects of
$ {\Cal M}(X/k,A) $ contain
$ l $-adic sheaves as in (1.1), then an objects of
$ {\Cal M}(A)_{\langle k \rangle} $ contains the inductive limit of
$ l $-adic sheaves which has a Galois action in the weak sense.
(This is however enough to define the weight filtration.)
Indeed, the
$ l $-adic part of an object of
$ {\Cal M}(K,A)_{\langle k \rangle} $ is identified with a
$ {\Bbb Q}_{l} $-vector space having an action of
$ \Gal(\overline{K}/K) $ (unramified over
$ \Spec R $ for some finitely generated smooth
$ k $-subalgebra
$ R $ of
$ K $ whose fraction field is
$ K) $.
Furthermore the pull-back by
$ K' \supset K $ is given by the composition of natural morphisms
$ \Gal(\overline{K'}/K') \rightarrow \Gal(\overline{K}K'/K') =
\Gal(\overline{K}
/\overline{K}\cap K') \rightarrow \Gal(\overline{K}/K) $.

(iv) Let
$ \MHM(X,A)^{\go} $,
$ \MHS(A)^{\go} $ be the categories of mixed Hodge Modules or structures of
geometric origin.
See [40].
We can define similarly
$ \MHM(X,A)_{\langle k \rangle }^{\go} $,
$ \MHS(A)_{\langle k \rangle }^{\go} $ as full subcategories of
$ \MHM(X,A)_{\langle k \rangle} $,
$ \MHS(A)_{\langle k \rangle} $ using morphisms of complex varieties and
the direct image and pull-back functors between
$ \MHM(X,A)_{\langle k \rangle} $.
Then we have a natural functor
$$
\MHM(X,A)_{\langle k \rangle }^{\go} \rightarrow \MHM(X,A)^{\go},
\quad
\MHS(A)_{\langle k \rangle }^{\go} \rightarrow \MHS(A)^{\go}.
\leqno(2.5.1)
$$
These are essentially surjective.
But it is not clear whether they are fully faithful.
It would be true for morphisms between pure objects, if the Hodge
conjecture is true.
But the mixed case is not clear.

In (2.1) we can restrict to the subcategory consisting of objects of
geometric origin, because
$ {\Cal M}(X,A)_{\langle k \rangle }^{\go} =
\hbox{\rm ind\,lim} \,{\Cal M}(X_{R'}/k_{R'},A)^{\go}[-d_{R'}] $.
In this case the surjectivity of the cycle map is equivalent to the Hodge
type conjecture.
See [37], [40].

\proclaim{{\bf 2.6.~Proposition}}
Let
$ K, K' $ be as in (2.1) such that
$ K \subset K' $,
$ \overline{k} \cap K = \overline{k} \cap K' $ and
$ K $,
$ K'$ are finitely generated over
$ k $.
Let
$ \pi ' : \Spec K' \rightarrow \Spec K $ denote the canonical morphism.
Then for a pure object
$ M $ in
$ \MHS(K,A)_{\langle k \rangle} $, the canonical morphism
$$
\Ext^{i}(A_{K,\langle k \rangle},M) \rightarrow \Ext^{i}(A_{K',\langle k
\rangle},
\pi ^{\prime *}M)
\leqno(2.6.1)
$$
is injective.
\endproclaim

\demo\nofrills {Proof.\usualspace}
$ Let $
$ k' = \overline{k} \cap K $.
Let
$ u $ be an element of the source, and assume its image in the target is
zero.
Then there exist finitely generated smooth
$ k' $-subalgebras
$ R \subset R' $ of
$ {\Bbb C} $ such that
$ K, K' $ are the fraction fields of
$ R, R' $, and
$ u $ is defined over
$ R $.
We may assume that there is a multivalued section of
$ S' := \Spec R' \rightarrow S := \Spec R $ such that its image
$ S'' $ is finite etale over
$ S $.
Let
$ \pi '' : S'' \rightarrow S $ denote the canonical morphism.
Then the image of
$ u $ in
$$
\Ext^{i}(A_{S''/k'},\pi ^{\prime \prime *}M) = \Ext^{i}(A_{S/k'},\pi
''_{*}\pi ^{\prime \prime *}M)
$$
is zero, and the injection
$ M \rightarrow \pi ''_{*}\pi ^{\prime \prime *}M $ splits by the
semisimplicity of pure objects.
So
$ u $ is zero, and the assertion follows.
\enddemo

\proclaim{{\bf 2.7.~Proposition}}
If
$ K $ is algebraically closed, we have equivalences of categories
$$
\MHM(X/K,A)_{\langle k \rangle} \rightarrow
\MHM(X/K,A)_{\langle \overline{k} \rangle},
\quad \MHS(K,A)_{\langle k \rangle} \rightarrow
\MHS(K,A)_{\langle \overline{k}\rangle}.
$$
induced by
$ \MHM(X_{R'}/k_{R'},A) \rightarrow \MHM(X_{R'}\otimes
_{k_{R'}}\overline{k}
/\overline{k},A) $.
\endproclaim

\demo\nofrills {Proof.\usualspace}
It is enough to verify that a filtered
$ {\Cal D} $-Module on
$ X_{R'}\otimes _{k_{R'}}\overline{k}/\overline{k} $ is defined on
$ X_{R'}\otimes _{k_{R'}}k'/k' $ for a sufficiently large algebraic
extension
$ k' $ of
$ k $ in
$ \overline{k} $ (and a similar assertion holds for morphisms).
\enddemo

\proclaim{{\bf 2.8.~Theorem}}
The category
$ {\Cal M}(X/K,A)_{\langle k \rangle} $ is an
$ (A \cap \overline{k}) $-linear abelian category (where
$ A = {\Bbb Q} $ unless
$ {\Cal M}(X/K,A)_{\langle k \rangle} = \MHM(X/K,A)_{\langle k \rangle}) $,
and the bounded derived categories
$ D^{b}{\Cal M}(X/K,A)_{\langle k \rangle} $ are stable by standard functors
like
$ f_{*}, f_{!}, f^{*}, f^{!}, $ etc. in a compatible way with the functor
$ \iota $.
\endproclaim

\demo\nofrills {Proof.\usualspace}
This follows from (1.2) using an analogue of the generic base change
theorem in [13].
\enddemo

\noindent
{\it Remark.} Any object of
$ {\Cal M}(X/K,A)_{\langle k \rangle} $ is represented by an object of
$ {\Cal M}(X_{R'/k_{R'}},A)[-d_{R'}] $ whose underlying perverse sheaf has
a Whitney stratification of
$ \overline{X}_{R'/k_{R'}} $ such that the strata are smooth over
$ R' $, where
$ \overline{X}_{R'/k_{R'}} $ is a smooth compactification of
$ X_{R'/k_{R'}} $ over
$ R' $.

\bigskip\bigskip
\centerline{{\bf 3.~Cycle Maps}}

\bigskip
\noindent
{\bf 3.1.}
Let
$ X $,
$ K $ be as in (2.1) with structure morphism
$ a_{X/K} : X \rightarrow \Spec K $.
We define
$$
\aligned
A_{X/K,\langle k \rangle}(j)
&= {a}_{X/K}^{*}A_{K,\langle k \rangle}(j) \in D^{b}{\Cal M}
(X/K,A)_{\langle k \rangle},
\\
H^{i}(X/K,A_{\langle k \rangle}(j))
&= H^{i}(a_{X/K})_{*}A_{X/K,\langle k \rangle}(j)
\in {\Cal M}(K,A)_{\langle k \rangle}.
\endaligned
$$
These are denoted by
$ A_{X,\langle k \rangle}(j) $,
$ H^{i}(X,A_{\langle k \rangle}(j)) $ to simplify the notation if
$ k = {\Bbb C} $.
This definition is compatible with the functor
$ \iota $, and
$ H^{i}(X/K,A_{\langle k \rangle}(j)) $ is the limit of
$ H^{i}\pi _{*}A_{X_{R}/k_{R}}(j) $, where
$ \pi : X_{R} \rightarrow S = \Spec R $ denotes the structure morphism.

If
$ A = {\Bbb Q} $ and
$ X $ is smooth over
$ K = {\Bbb C} $, let
$$
\Hdg^{p}(X,{\Bbb Q}_{\langle k \rangle})
= \Hom({\Bbb Q}_{\langle k \rangle},
H^{2p}(X, {\Bbb Q}_{\langle k \rangle}(p))).
$$
If the mixed sheaves in (1.1) are mixed Hodge Modules, this is called the
group of
$ k $-finite (or
$ k $-horizontal) Hodge cycles of codimension
$ p $.

\medskip
\noindent
{\bf 3.2.}
{\it Remark.} (i)
Assume the mixed sheaves in (1.1) are mixed Hodge Modules.
By Remark (i) of (2.5),
$ \Hdg^{p}(X,{\Bbb Q}_{\langle k \rangle}) $ is identified with the space
of
$ {\Bbb Q} $-Hodge cycles (in the usual sense) which are
$ k $-horizontal (i.e. annihilated by the connection over
$ \overline{k} $).
See also [22], [23].

When
$ A = {\Bbb Q} $, the Hodge conjecture implies that every
$ {\Bbb Q} $-Hodge cycle in the usual sense is a
$ \overline{\Bbb Q} $-finite (i.e.
$ \overline{\Bbb Q} $-horizontal) Hodge cycle.
It is expected that Hodge cycles are
$ \overline{\Bbb Q} $-horizontal.

\medskip
\noindent
{\bf 3.3.}
{\it Deligne cohomology.} Let
$ X $ be a proper smooth variety over
$ {\Bbb C} $, and
$ {H}_{\Cal D}^{i}(X,{\Bbb Q}(j)) $ the
$ {\Bbb Q} $-Deligne cohomology in the usual sense [17].
By [3] it is isomorphic to
$$
\Ext^{i}({\Bbb Q}^{H}, (a_{X})_{*}{\Bbb Q}_{X}^{H}(j)) =
\Ext^{i}({\Bbb Q}_{X}^{H}, {\Bbb Q}_{X}^{H}(j)),
$$
where
$ {\Bbb Q}^{H} $ and
$ {\Bbb Q}_{X}^{H} $ belong to
$ D^{b}\MHS({\Bbb Q}) $ and
$ D^{b}\MHM(X,{\Bbb Q}) $ respectively (see [36]), and the last isomorphism
follows from the adjoint relation between the direct image and the
pull-back.
We have by [11]
$$
J^{p}(X)_{\Bbb Q} = \Ext^{1}({\Bbb Q}^{H}, H^{2p-1}(X, {\Bbb Q}
(p))),
\leqno(3.3.1)
$$
where
$ J^{p}(X)_{\Bbb Q} $ is Griffiths' intermediate Jacobian [24] tensored
with
$ {\Bbb Q} $.

We define
$ {H}_{\Cal D}^{i}(X,{\Bbb Q}_{\langle k \rangle}(j)) $ to be
$$
\Ext^{i}({\Bbb Q}_{\langle k \rangle}, (a_{X})_{*}{\Bbb Q}_{X,\langle k
\rangle}(j))
=\Ext^{i}({\Bbb Q}_{X,\langle k \rangle}, {\Bbb Q}_{X,\langle k
\rangle}(j)),
$$
where the extension group is taken in the derived category of
$ {\Cal M}(A)_{\langle k \rangle} $ or
$ {\Cal M}(X,A)_{\langle k \rangle} $.
This is the limit of
$ {H}_{\Cal D}^{i}(X_{R}/k_{R},{\Bbb Q}_{k}(j)) $ which is defined to be
$$
\Ext^{i}({\Bbb Q}_{k_{R}}, (a_{X_{R}/k_{R}})_{*}{\Bbb Q}_{X_{R}/k_{R}}(j))
=\Ext^{i}({\Bbb Q}_{X_{R}/k_{R}}, {\Bbb Q}_{X_{R}/k_{R}}(j)).
\leqno(3.3.2)
$$
See (1.3).
By (2.1.4) we have a natural morphism
$ {H}_{\Cal D}^{i}(X,{\Bbb Q}_{\langle k \rangle}(j)) \rightarrow
{H}_{\Cal D}^{i}(X,{\Bbb Q}(j)) $.
Its image is denoted by
$ {H}_{\Cal D}^{i}(X,{\Bbb Q}(j))^{\langle k \rangle} $, and is called the
$ k $-finite part of the Deligne cohomology associated with the mixed
sheaves
$ {\Cal M} $.
Since we have a noncanonical isomorphism
$$
(a_{X})_{*}{a}_{X}^{*}{\Bbb Q}_{\langle k \rangle}
\simeq \oplus _{j}H^{j}(X,{\Bbb Q}_{\langle k \rangle})[-j]
$$
in the derived category of
$ {\Cal M}(A)_{\langle k \rangle} $ (see [37], (8.1.7)), there exists a
decreasing filtration
$ L $ on
$ {H}_{\Cal D}^{i}(X,{\Bbb Q}_{\langle k \rangle}(j)) $ induced by the
truncation
$ \tau $ on
$ (a_{X})_{*}{\Bbb Q}_{X,\langle k \rangle}(j) $ such that
$$
{\Gr}_{L}^{r}{H}_{\Cal D}^{i}(X,{\Bbb Q}_{\langle k \rangle}(j)) =
\Ext^{r}({\Bbb Q}_{\langle k \rangle}, H^{i-r}(X,{\Bbb Q}_{\langle k
\rangle}(j))).
\leqno(3.3.3)
$$
See [37], (8.1).
For
$ (i,j) = (2p,p) $, we have a morphism of short exact sequences
$$
\CD
0 @>>> J^{p}(X)_{\Bbb Q}^{\langle k \rangle } @>>>
{H}_{{\Cal D}}^{2p}(X,{\Bbb Q}(p))^{\langle k \rangle} @>>>
\Hdg^{p}(X,{\Bbb Q}_{\langle k \rangle}) @>>> 0
\\
@. @VVV @VVV @VVV @.
\\
0 @>>> J^{p}(X)_{\Bbb Q} @>>> {H}_{{\Cal D}}^{2p}(X,{\Bbb Q}(p)) @>>>
\Hdg^{p}(X,{\Bbb Q}) @>>> 0
\endCD
\leqno(3.3.4)
$$
where
$ J^{p}(X)_{\Bbb Q}^{\langle k \rangle } $ is defined to be the image of
$$
\Ext^{1}({\Bbb Q}_{\langle k \rangle}, H^{2p-1}(X, {\Bbb Q}_{\langle k
\rangle}(p)))\rightarrow \Ext^{1}({\Bbb Q}^{H}, H^{2p-1}(X, {\Bbb Q}(p)))
= J^{p}(X)_{\Bbb Q}.
\leqno(3.3.5)
$$
{\it Remark.} It is not clear whether (3.3.5) is injective.
See Remark (ii) of (2.5).

\medskip
\noindent
{\bf 3.4.}
{\it Cycle map.} For a smooth proper complex variety
$ X $, let
$ \CH^{p}(X)_{\Bbb Q} $ be the Chow group of codimension
$ p $ cycles on
$ X $ tensored with
$ {\Bbb Q} $.
By [36], [37] we have a cycle map
$$
cl : \CH^{p}(X)_{\Bbb Q} \rightarrow {H}_{\Cal D}^{2p}(X,{\Bbb
Q}_{\langle k \rangle}(p)).
\leqno(3.4.1)
$$
This is the limit of the cycle map of (1.3.2) :
$$
cl_{R} : \CH^{p}(X_{R})_{\Bbb Q} \rightarrow {H}_{\Cal D}^{2p}
(X_{R}/k_{R},{\Bbb Q}_{k}(p)).
\leqno(3.4.2)
$$

We have the induced filtration
$ L $ (see (3.3.3)) on
$ \CH^{p}(X)_{\Bbb Q} $, and the graded pieces of the cycle map
$$
{\Gr}_{L}^{r}cl : {\Gr}_{L}^{r}\CH^{p}(X)_{\Bbb Q} \rightarrow
\Ext^{r}({\Bbb Q}
_{\langle k \rangle}, H^{2p-r}(X,{\Bbb Q}_{\langle k \rangle}(p))
\leqno(3.4.3)
$$
is called the higher Abel-Jacobi map.
It may be expected that
$ L $ gives a conjectural filtration of Beilinson and Bloch (see [4], [7]
and also [33]) on the Chow group.
See for example [1].

For a
$ k $-variety
$ X $ as in (2.1), we have similarly the induced filtration
$ L $ on
$ \CH^{p}(X)_{\Bbb Q} $ and
$$
{\Gr}_{L}^{r}cl_{K} : {\Gr}_{L}^{r}\CH^{p}(X)_{\Bbb Q} \rightarrow
\Ext^{r}({\Bbb Q}
_{K,\langle k \rangle}, H^{2p-r}(X/K,{\Bbb Q}_{\langle k \rangle}(p)).
\leqno(3.4.4)
$$
This gives (3.4.3) by passing to the limit.

We have a short exact sequence
$$
0 \rightarrow L^{1}\CH^{p}(X)_{\Bbb Q} \rightarrow
\CH^{p}(X)_{\Bbb Q} \rightarrow {\Gr}_{L}^{1}\CH^{p}(X)_{\Bbb Q}
\rightarrow 0,
\leqno(3.4.5)
$$
and the cycle map induces a morphism of this exact sequence to the first
row of the commutative diagram (3.3.4).
We can verify that the morphism
$ L^{1}\CH^{p}(X)_{\Bbb Q} \rightarrow J^{p}(X)_{\Bbb Q} $ coincides
with Griffiths' Abel-Jacobi map (by using [17]).
See [36], (4.5.20).
Indeed, we can prove that, for
$ \zeta \in \CH^{p}(X)_{\Bbb Q} $ represented by a closed irreducible
subvariety
$ Z $ of
$ X $,
$ cl(\zeta ) $ is represented by the composition
$$
{\Bbb Q}^{H} \rightarrow {H}_{Z}^{2p}(X,{\Bbb Q})(p) \rightarrow
(a_{X})_{*}{\Bbb Q}_{X}^{H}(p)[2p].
$$
(The last has been known to specialists as Deligne's cycle map.)

\medskip
\noindent
{\it Remark.} We can define a cycle map from Bloch's higher Chow group
$ \CH^{p}(X,n)_{\Bbb Q} $ [8] to
$ \Ext^{2p-n}({\Bbb Q}_{X,\langle k \rangle}, {\Bbb Q}_{X,\langle k
\rangle}(p)) $, and get the filtration
$ L $.
See [37], (8.3).

\medskip
\noindent
{\bf 3.5.}
{\it Correspondence.} For smooth proper
$ K $-varieties
$ X, Y $, we define the group of correspondences by
$$
C^{i}(X,Y)_{\Bbb Q} = \CH^{i+\dim X}(X\times_{K} Y)_{{\Bbb Q},}
$$
where
$ X $ is assumed to be equidimensional.
Then we have the cycle map
$$
\aligned
C^{i}(X,Y)_{\Bbb Q}
&\rightarrow \Ext^{2i+2\dim X}({\Bbb Q}_{X\times Y/K,\langle k \rangle},
{\Bbb Q}_{X\times Y/K,\langle k \rangle}(i+\dim X))
\\
&= \Hom((a_{X/K})_{*}{\Bbb Q}_{X/K,\langle k \rangle}, (a_{Y/K})_{*}
{\Bbb Q}_{Y/K,\langle k \rangle}(i)[2i]).
\endaligned
\leqno(3.5.1)
$$
This cycle map is compatible with the composition of correspondences.
See [40,~II], (3.3).
This implies that the action of
$ C^{i}(X,Y)_{\Bbb Q} $ on the Chow groups
corresponds by the cycle map (3.4.1) to the composition of morphisms with
the image of (3.5.1), because
$ C^{i}(\Spec K,X)_{\Bbb Q} = \CH^{i}(X)_{\Bbb Q} $.
So we get a commutative diagram
$$
\CD
\CH^{p}(X)_{\Bbb Q} @>{cl}>>
\Hom({\Bbb Q}_{K,\langle k \rangle}, (a_{X/K})_{*}
{\Bbb Q}_{X/K,\langle k \rangle}(p)[2p])
\\
@VV{\Gamma_{*}}V @VV{\Gamma_{*}}V
\\
\CH^{p+i}(Y)_{\Bbb Q}
@>{cl}>>
\Hom({\Bbb Q}_{K,\langle k \rangle}, (a_{Y/K})_{*}
{\Bbb Q}_{Y/K,\langle k \rangle}(p+i)[2p+2i])
\endCD
\leqno(3.5.2)
$$
for a correspondence
$ \Gamma \in C^{i}(X,Y)_{\Bbb Q} $.
Since
$$
\Gamma _{*} : (a_{X/K})_{*}{\Bbb Q}_{X/K,\langle k \rangle} \rightarrow
(a_{Y/K})_{*}{\Bbb Q}_{Y/K,\langle k \rangle} (i)[2i]
$$
preserves the filtration
$ \tau $, we see that the action of a correspondence on the Chow groups
preserves the filtration
$ L $, and the action on the graded pieces depends only on the cohomology
class of the correspondence.

\medskip
\noindent
{\it Remark.} If
$ \zeta \in C^{i}(X,Y)_{\Bbb Q} $ is represented by an irreducible
closed subvariety
$ Z $, let
$ Z' \rightarrow Z $ be a resolution of singularities, and
$ p_{i} $ denotes the composition with the projection to each direct factor.
Then by [loc.~cit] the image of
$ \zeta $ by (3.5.1) is the composition of the restriction and Gysin
morphisms :
$$
\aligned
{p}_{1}^{*} : (a_{X/K})_{*}{\Bbb Q}_{X/K,\langle k \rangle}
&\rightarrow (a_{Z'/K})_{*}{\Bbb Q}_{Z'/K,\langle k \rangle},
\\
p_{2*} : (a_{Z'/K})_{*}{\Bbb Q}_{Z'/K,\langle k \rangle}
&\rightarrow (a_{Y/K})_{*}{\Bbb Q}_{Y/K,\langle k \rangle} (i)[2i].
\endaligned
$$

\proclaim{{\bf 3.6.~Proposition}}
Let
$ X $ be a smooth proper
$ K $-variety of dimension
$ n $, where
$ K $ is as in (2.1).
Then
$ L^{2}\CH^{n}(X)_{\Bbb Q} $ (see (3.4.4)) coincides with the kernel
of the Albanese map, and we get an injection
$$
\Alb_{X/K}(K) \rightarrow \Ext^{1}({\Bbb Q}_{K,\langle k \rangle},
H^{2n-1}(X/K,{\Bbb Q}_{\langle k \rangle}(p)).
$$
\endproclaim

\demo\nofrills {Proof.\usualspace}
Since
$ \CH^{n}(X)_{\Bbb Q}, \Alb_{X/K}(K) $, and
$ H^{2n-1}(X/K,{\Bbb Q}_{\langle k \rangle}(p)) $ do not change by a
birational morphism of smooth proper varieties (see [19], [29]), the
assertion is reduced to the projective case.

It is clear that
$ L^{2}\CH^{n}(X)_{\Bbb Q} $ is contained in the kernel of the
Albanese map.
(Indeed, restricting to the fiber over the point corresponding to the
inclusion
$ K \rightarrow {\Bbb C} $ and using (1.1.2), the assertion is reduced to
the case
$ K = {\Bbb C} $.)
Let
$ \zeta \in L^{1}\CH^{n}(X)_{\Bbb Q} $ belonging to the kernel of the
Albanese map.
Let
$ \Gamma $ be the projector in
$ \CH^{n}(X\times _{K}X)_{\Bbb Q} $ corresponding to the Albanese
motive in [32], 4.1 (i.e.
$ \Gamma ^{2} = \Gamma $ and the action on
$ H^{i}(X\otimes _{K}\overline{K},{\Bbb Q}_{l}) $ is the identity for
$ i = 2n - 1 $, and zero otherwise).
Furthermore
$ \Gamma _{*}\zeta = 0 $ by loc.~cit.
We have to show the vanishing of
$$
cl(\zeta ) \in \Ext^{1}({\Bbb Q}_{K,\langle k \rangle},H^{2n-1}(X/K,{\Bbb
Q}
_{\langle k \rangle})).
$$
But this follows from (3.5.2) applied to
$ \Gamma $ and
$ \zeta $, because the action of
$ \Gamma _{*} $ on
$ H^{2n-1}(X/K,{\Bbb Q}_{\langle k \rangle}) $ is the identity.
So we get the assertion.
\enddemo

\noindent
{\bf 3.7.} {\it Example.} Let
$ X $ be a smooth subvariety in
$ {\Bbb P}_{{\Bbb C}}^{3} $ which is a geometric generic fiber of a pencil
$ f : Y_{k} \rightarrow {\Bbb P}_{k}^{1} $.
Here
$ Y_{k} $ is the blowing-up of
$ {\Bbb P}_{k}^{3} $ along a smooth curve
$ C_{k} $ which is the intersection of two smooth surfaces of degree
$ d $,
and the fibers of
$ f $ are the members of the linear system defined by the two surfaces
(i.e. it is a Lefschetz pencil by embedding
$ {\Bbb P}_{k}^{3} $ into a projective space using the ample line bundle
$ {\Cal O}_{{\Bbb P}^{3}}(d) $).
We assume that each (geometric) singular fiber has at most one singular
point and it is an ordinary double point.
We choose an embedding of
$ K := k({\Bbb P}_{k}^{1}) $ into
$ {\Bbb C} $ so that the base change of the generic fiber is
$ X $.

Let
$ Y, C $ denote the base change of
$ Y_{k}, C_{k} $ by
$ k \rightarrow {\Bbb C} $.
Then
$$
H^{j}(Y,{\Bbb Q}) = H^{j}({\Bbb P}^{3},{\Bbb Q}) \oplus
H^{j-2}(C,{\Bbb Q})(-1).
\leqno(3.7.1)
$$
Let
$ {\CH}_{\hom}^{p}(Y)_{{\Bbb Q}} $ denote the subgroup of
$ \CH^{p}(Y)_{{\Bbb Q}} $ consisting of homologically equivalent to zero
cycles.
Then we have a commutative diagram
$$
\CD
{\CH}_{\hom}^{1}(C)_{{\Bbb Q}} @>{\simeq}>>
{\CH}_{\hom}^{1}(C\times {\Bbb P}^{1})_{{\Bbb Q}} @>{\simeq}>>
{\CH}_{\hom}^{2}(Y)_{{\Bbb Q}}
\\
@VV{\simeq}V @VV{\simeq}V @VV{\simeq}V
\\
J^{1}(C)_{{\Bbb Q}} @>{\simeq}>>
J^{1}(C\times {\Bbb P}^{1})_{{\Bbb Q}} @>{\simeq}>>
J^{2}(Y)_{{\Bbb Q}}
\endCD
\leqno(3.7.2)
$$
where the horizontal morphisms are the pull-back by the projection
$ C\times {\Bbb P}^{1} \rightarrow C $ and the direct image by the inclusion
$ C\times {\Bbb P}^{1} \rightarrow Y $, and all the morphisms are
isomorphisms.
This implies the isomorphisms
$$
{\CH}_{\hom}^{1}(C_{k})_{{\Bbb Q}} \simto
{\CH}_{\hom}^{1}(C_{k}\times _{k}{\Bbb P}_{k}^{1})_{{\Bbb Q}} \simto
{\CH}_{\hom}^{2}(Y_{k})_{{\Bbb Q}}
\leqno(3.7.3)
$$
together with the injective morphisms to the top row of the above diagram.
In particular, we get the injectivity of the Abel-Jacobi map
$$
\CH_{\hom}^{2}(f^{-1}(U_{k}))_{{\Bbb Q}} \rightarrow
J^{2}(f^{-1}(U))_{{\Bbb Q}} :=\Ext_{\MHS}^{1}({\Bbb Q},
H^{3}(f^{-1}(U),{\Bbb Q})(2))
\leqno(3.7.4)
$$
for
$ U_{k} = {\Bbb P}_{k}^{1} $.
We see that this holds for any open subvarieties
$ U_{k} $ of
$ {\Bbb P}_{k}^{1} $ using the diagram
$$
\CD
{\CH}_{\hom}^{1}(Y_{\Sigma})_{{\Bbb Q}}^{(0)} @>>>
{\CH}_{\hom}^{1}(Y)_{{\Bbb Q}} @>>>
{\CH}_{\hom}^{2}(Y_{U})_{{\Bbb Q}}
\\
@VVV @VVV @VVV
\\
\Hom_{\MHS}({\Bbb Q},H^{2}(Y_{\Sigma},{\Bbb Q})^{(0)}(1)) @>>>
J^{2}(Y)_{{\Bbb Q}} @>>>
J^{2}(Y_{U})_{{\Bbb Q}}
\endCD
$$
where
$ Y_{U} = f^{-1}(U) $,
$ Y_{\Sigma} = f^{-1}(\Sigma) $ with
$ \Sigma = {\Bbb P}^{1} \setminus U $, and the superscript
$ {}^{(0)} $ means the kernel of the morphism to
$ H^{4}(Y,{\Bbb Q})(2) $.
(Indeed, the left vertical morphism is surjective
by the Lefschetz theorem applied to the blowing-up of
$ Y_{s} := f^{-1}(s) $ at singular points).
Note that the fibers are rational homology manifolds
so that the cohomology, intersection cohomology and homology of
$ Y_{s} $ all coincide in this case.
Here it is also possible to replace
$ Y $ with its blowing-up at the singular points of fibers, and argue
as in (4.7) (but the argument is easier, because
$ H^{1}(Y_{s},{\Bbb Q}) = 0 $.)

So we can control codimension two cycles on
$ X $ which are defined over
$ K $.
But it is not easy to extend this to all the cycles even in this simple
example, because we have to consider arbitrary extensions of
$ K $.
We will reduce the problem to the injectivity of certain Abel-Jacobi map
for varieties over number fields.

As for the relation with the filtration
$ L $ in (3.4), we have
$ R^{j}f_{*}{\Bbb Q}_{Y} = 0 $ for
$ j $ odd, and
$ R^{2}f_{*}{\Bbb Q}_{Y} $ is an intersection complex with local system
coefficients (i.e. the intermediate direct image [6]) which underlies
naturally a pure Hodge Module.
So we have an isomorphism of mixed Hodge structures
$$
H^{3}(f^{-1}(U),{\Bbb Q}) = H^{1}(U,R^{2}f_{*}{\Bbb Q}_{Y}),
\leqno(3.7.5)
$$
and this implies
$$
\Ext^{2}({\Bbb Q}_{U},R^{2}f_{*}{\Bbb Q}_{Y}) =
J^{2}(f^{-1}(U))_{{\Bbb Q}}.
\leqno(3.7.6)
$$
Here we assume that the mixed sheaves in (1.1) are mixed Hodge Modules.
Then (3.7.4) gives part of of
$ {\Gr}_{L}^{2}cl $ (because
$ {\Gr}_{L}^{1}cl = 0 $ by
$ H^{3}(X,{\Bbb Q}) = 0 $), and we get the nonvanishing of
$ {\Gr}_{L}^{2}cl $ in this case by (2.6).
Note that the nonvanishing follows also from Bloch's diagonal argument
easily in the case
$ H^{3}(X,{\Bbb Q}) = 0 $ as remarked after (0.2) in the introduction.

\bigskip\bigskip
\centerline{{\bf 4.~Proof of Main Theorems}}

\bigskip
\noindent
In this section we prove the main theorems for mixed sheaves in (1.1).

\medskip
\noindent
{\bf 4.1.}
{\it Proof of} (0.3).
It is enough to show the surjectivity of
$$
\CH^{p}(X)_{\Bbb Q} \rightarrow
\Hdg^{p}(X,{\Bbb Q}_{\langle k \rangle}), \quad L^{1}\CH^{p}(X)_{\Bbb Q}
\rightarrow J^{p}(X)_{\Bbb Q}^{\langle k \rangle}.
$$
Let
$ R $ be a finitely generated smooth
$ k $-subalgebra of
$ {\Bbb C} $ with a smooth proper morphism
$ \pi : X_{R} \rightarrow S := \Spec R $ such that
$ X_{R}\otimes _{R}{\Bbb C} = X $.
Put
$ A = {\Bbb Q} $.
Then
$ H^{i}(X,{\Bbb Q}_{\langle k \rangle}(j)) $ is represented by a variation
of Hodge structure
$ H^{i}\pi _{*}A_{X_{R}/k_{R}}(j) $ on
$ S $.
For
$ i = 0, 1 $, we have
$$
\Ext^{i}(A_{S/k_{R}}, H^{2p-i}\pi _{*}A_{X_{R}/k_{R}}(p)) = \Hom(A_{k_{R}},
H^{i}(a_{S/k_{R}})_{*}H^{2p-i}\pi _{*}A_{X_{R}/k_{R}}(p)).
$$
by using the Leray spectral sequence
$$
{E}_{2}^{i,j} = \Ext^{i}(A_{k_{R}}, H^{j}(a_{S/k_{R}})_{*}M) \Rightarrow
\Ext^{i+j}(A _{S/k_{R}}, M),
\leqno(4.1.1)
$$
where
$ M = H^{2p}\pi _{*}A_{X_{R}/k_{R}}(p) $ or
$ H^{2p-1}\pi _{*}A_{X_{R}/k_{R}}(p) $ (because
$ {E}_{2}^{i,0} = 0 $ for
$ i = 1, 2 $ in the latter case by assumption).
Since we have a noncanonical isomorphism (see [37], (6.10)) :
$$
\pi _{*}A_{X_{R}/k_{R}}
\simeq \oplus _{i} (H^{i}\pi _{*}A_{X_{R}/k_{R}})[-i]\quad
\text{in }D^{b}\MHM(S/k_{R},A),
\leqno(4.1.2)
$$
we can identify
$ H^{i}(a_{S/k_{R}})_{*}H^{2p-i}\pi _{*}A_{X_{R}/k_{R}}(p) $ noncanonically
with a subobject (and canonically with a subquotient) of
$ H^{2p}(X_{R}/k_{R}, A(p)) $.
So the assertion is reduced to the Hodge conjecture for a smooth
compactification of
$ X_{R}\otimes _{k_{R}}{\Bbb C} $.
See Remark (i) below.
Then it is further reduced to the Hodge conjecture for a smooth projective
variety due to Chow's lemma.

\medskip
\noindent
{\it Remarks.} (i) Let
$ X_{k} $ be a smooth
$ k $-variety such that
$ k = \overline{k} \cap k(X_{k}) $ (i.e.
$ X_{k}/k $ is geometrically irreducible).
Let
$ \overline{X}_{k} $ be a smooth compactification of
$ X_{k} $ such that
$ \overline{X}_{k} \setminus X_{k} $ is a divisor with normal crossings.
Then the cycle map
$$
\CH^{p}(X)_{\Bbb Q} \rightarrow \Hom({\Bbb Q}_{k},H^{2p}(X _{k}/k,{\Bbb
Q}_{k}(p)))
\leqno(4.1.3)
$$
is surjective if the Hodge conjecture for codimension
$ p $ cycles on
$ \overline{X}_{\Bbb C} := \overline{X}_{k}\otimes _{k}{\Bbb C} $ is
true.

Indeed,
$ H^{2p}(X_{k}/k,{\Bbb Q}_{k}) $ has weights
$ \ge 2p $ and
$ {\Gr}_{2p}^{W}H^{2p}(X_{k}/k,{\Bbb Q}_{k}) $ is a quotient of
$ H^{2p}(\overline{X}_{k}/k,{\Bbb Q}_{k}) $.
See [12].
So the assertion is reduced to the case
$ X_{k} $ is smooth proper.
Then the assertion is more or less well-known (using the Galois action on
$ l $-adic cohomology or de Rham cohomology).
See for example [38], (8.5,iii).

(ii) For a smooth proper
$ k $-variety
$ X_{k} $ and a finitely generated smooth
$ k $-subalgebra
$ R $ of
$ {\Bbb C} $, let
$ X = X_{k}\otimes _{k}{\Bbb C}, X_{R} = X_{k}\otimes _{k}R $, and
$ S = \Spec R $.
Assume
$ k = \overline{k} \cap k(X_{k}) $ and the mixed sheaves in (1.1) are
mixed Hodge Modules.
Then
$ J^{p}(X)^{\langle k \rangle} = J^{p}(X) $, and the Abel-Jacobi map to
$ J^{p}(X)^{\langle k \rangle} $ is not surjective if
$ H^{2p-1}(X) $ has level
$ > 1 $.

\medskip
\noindent
{\bf 4.2.}
{\it Proof of} (0.1).
Let
$ R $ be a finitely generated smooth
$ k $-subalgebra of
$ {\Bbb C} $ with a smooth proper morphism
$ \pi : X_{R} \rightarrow S := \Spec R $ such that
$ X_{R}\otimes _{R}{\Bbb C} = X $.
We may assume
$ X $ is connected.
Put
$$
M = H^{2n-j}\pi _{*}A_{X_{R}/k_{R}}(n).
$$
Let
$ R' $ be the affine ring of an affine open subvariety
$ U $ of
$ X_{R} $ (i.e.
$ U = \Spec R') $.
Put
$ S' = \Spec R' $ with
$ \pi ' : S' \rightarrow S $ the canonical morphism.
Let
$ R^{(m)} $ be the
$ m $-ple tensor product of
$ R' $ over
$ R $.
We choose an embedding
$ R^{(m)} \rightarrow {\Bbb C} $ over
$ R $.
(Here
$ R^{(m)}\otimes _{k_{R}}\overline{k} $ is integral so that
$ R^{(m)} \cap \overline{k} = k_{R} $, because
$ \pi ' $ is smooth, and the generic fiber is geometrically irreducible.)
Put
$ S^{(m)} = \Spec R^{(m)} $ with
$ \pi ^{(m)} : S^{(m)} \rightarrow S $ the canonical morphism.
For a nonempty affine open subvariety
$ S'' $ of
$ S^{(m)} $, let
$ \pi '' : S'' \rightarrow S $ denote the restriction of
$ \pi ^{(m)} $.
It is enough to show that
$$
\hbox{\rm ind\,lim}_{S''} \Ext^{j}(A_{S''/k_{R}}, \pi ^{\prime \prime *}M)
\leqno(4.2.1)
$$
with
$ S'' $ running over nonempty affine open subvarieties of
$ S^{(m)} $, has at least dimension
$ m $, and the morphism induced by the base change associated with the
first functor of (2.1.1) has at least rank
$ m $, because the effect by the second functor of (2.1.1) is covered by
(2.6).
(Concerning the first functor of (2.1.1), we can do the same construction
as above after taking the base change by a finite extension
$ k_{R'} \rightarrow k_{R''} $.)

By the edge morphism of the Leray spectral sequence, we have a canonical
morphism
$$
\Ext^{j}(A_{S''/k_{R}}, \pi ^{\prime \prime *}M) \rightarrow
\Hom(A_{S/k_{R}}, H^{j}\pi ''_{*}\pi ^{\prime \prime *}M)
$$
This is compatible with the transition morphism of the inductive limit, and
is surjective by the Lemma below.
So we may replace (4.2.1) with
$$
\hbox{\rm ind\,lim}_{S''} \Hom(A_{S/k_{R}}, H^{j}\pi ''_{*}\pi ^{\prime
\prime *}M),
\leqno(4.2.2)
$$
and we have to show that it has at least dimension
$ m $, and the morphism induced by the first functor of (2.1.1) has at
least rank
$ m $.

We first consider the case
$ m = 1 $, where we may assume
$ S'' = S' $ shrinking
$ U $ if necessary.
Let
$$
M' = H^{2n-j}\pi '_{!}A_{U/k_{R}}(n),\quad M^{\prime *} = H^{j}\pi
'_{*}A_{U/k_{R}},
$$
where
$ \pi ' $ is identified with the restriction of
$ \pi $ to
$ U $.
We may assume that the underlying
$ {\Bbb Q} $-complexes of
$ M' $ and
$ M^{\prime *} $ are local systems (shrinking
$ S $ if necessary).
Then they are dual of each other (up to a Tate twist and a shift of
complex), and
$$
H^{j}\pi '_{*}\pi ^{\prime *}M = M^{\prime *}\otimes M = {\Cal H}om (M',
M).
$$
We have a canonical element
$ \xi $ in
$$
\Hom(A_{S/k_{R}}, {\Cal H}om(M', M)) = \Hom(M', M),
$$
corresponding to the canonical morphism
$ M' \rightarrow M $.
We see that this is nonzero for any
$ U \subset X_{R} $ by restricting to the fiber at the geometric generic
point of
$ S $ and using the functor (2.1.4) together with the level of Hodge
structure, because
$ \Gamma (X, {\Omega }_{X}^{j}) \ne 0 $.
The compatibility with the base change is clear.

Now we consider the case
$ m > 1 $.
Since
$ {\pi }_{*}^{(m)}A_{S^{(m)}/k_{R}} $ is the
$ m $-ple tensor product of
$ \pi '_{*}A_{S'/k_{R}} $, we get a canonical injection
$$
\overset m \to\oplus (M^{\prime *}\otimes M) \rightarrow
H^{j}{\pi}_{*}^{(m)}A_{S^{(m)}/k_{R}}
\otimes M.
$$
We can verify that this gives an
$ m $-dimensional subspace of
$$
\Hom(A_{S/k_{R}}, H^{j}{\pi }_{*}^{(m)}A_{S^{(m)}/k_{R}}\otimes M),
$$
and of (4.2.2) using (2.1.4) and the level of Hodge structure as above.

\proclaim{{\bf 4.3.~Lemma}}
Let
$ X, S $ be smooth
$ k $-schemes of finite type with a smooth morphism
$ f : X \rightarrow S $ such that
$ k $ is algebraically closed in the function field of
$ S $ and
$ H^{j}f_{*}A_{X/k} $ is smooth (i.e. its underlying
$ {\Bbb Q} $-complex is a local system).
Let
$ M \in {\Cal M}(S/k, A)[-\dim S] $ with pure weight
$ - j $.
Then the canonical morphism
$$
\Ext^{j}(A_{X/k}, f^{*}M) = \Ext^{j}(A_{S/k}, f_{*}f^{*}M) \rightarrow
\Hom(A _{S/k}, H^{j}f_{*}f^{*}M)
$$
is surjective by shrinking
$ S $ if necessary.
\endproclaim

\demo\nofrills {Proof.\usualspace}
We may assume that
$ H^{j}f_{*}A_{X/k} $ and
$ M $ are smooth (i.e. their underlying
$ {\Bbb Q} $-complexes are local systems) by shrinking
$ S $ if necessary.
Let
$ \overline{f} : \overline{X} \rightarrow S $ be a smooth compactification
of
$ f $ (i.e.,
$ \overline{f} $ is smooth proper).
Then the canonical morphism
$ H^{j}\overline{f}_{*}\overline{f}^{*}M \rightarrow
{\Gr}_{0}^{W}H^{j}f_{*}f ^{*}M $ is surjective, and
$ W_{-1}H^{j}f_{*}f^{*}M = 0 $.
This implies the surjectivity of
$$
\Hom(A_{S/k}, H^{j}\overline{f}_{*}\overline{f}^{*}M) \rightarrow \Hom
(A_{S/k}, H^{j}f_{*}f^{*}M),
$$
and the assertion is reduced to the case where
$ f $ is smooth proper.
Then it follows from the decomposition (4.1.2) applied to the direct image
of
$ \overline{f}^{*}M $ by
$ \overline{f} $.
\enddemo

\medskip
\noindent
{\bf 4.4.}
{\it Proof of} (0.2).
We may assume
$ X $ is connected.
It is enough to construct cycles in
$ L^{2}\CH^{n}(X)_{\Bbb Q} $ whose images in (4.2.2) coincide with the
elements constructed in (4.2).
We first consider the case
$ X $ is projective.
With the notation of (4.2), let
$ K $ be the fraction field of
$ R $, and
$ X_{K} = X_{R}\otimes _{R}K $.
Let
$ K' $ be the function field of
$ X_{K} $, and
$ X_{K'} = X_{K}\otimes _{K}K' $.
By Murre's result [32], [33], there is a middle dimensional cycle
$ \zeta $ on
$ X_{K}\times _{K}X_{K} $ such that
$ \zeta \scirc \zeta = \zeta $ (as correspondences) and the action of its
base change
$ \zeta _{\Bbb C} \in \CH^{n}(X\times X) $ on
$ H^{j}(X,{\Bbb Q}) $ (via the correspondence) is the identity for
$ j \le 2n - 2 $, and zero otherwise.
Let
$ \zeta ' $ be the restriction of
$ \zeta $ to the generic fiber
$ X_{K'} $ of the first projection of
$ X_{K}\times _{K}X_{K} $ to
$ X_{K} $.
Let
$ \zeta _{K'} $ be the pull-back of
$ \zeta $ to
$ X_{K'}\times _{K'}X_{K'} $.
Then
$ (\zeta _{K'})_{*}\zeta ' = \zeta ' $ by
$ \zeta \scirc \zeta = \zeta $ (see Lemma below), and hence
$ \zeta ' \in L^{2}\CH^{n}(X_{K'})_{\Bbb Q} $ (using a diagram similar
to (3.5.2)).
Here we may assume that
$ \zeta $ is defined on
$ X_{R}\times _{R}X_{R} $ replacing
$ R $ if necessary.
We see that the image of
$ \zeta $ in
$ \Hom(A_{S/k}, H^{2}\pi ''_{*}\pi ^{\prime \prime *}M) $ in the notation
of (4.2) coincides with the canonical element constructed in (4.2), using
the action of
$ \zeta _{\Bbb C} $ on
$ H^{2n-2}(X,{\Bbb Q}) $ together with(2.1.4).
(A similar argument works for
$ L^{r}\CH^{n}(X)_{\Bbb Q} $ if the K\"unneth decomposition in the Chow
group [33] holds.)

To show that the image of the second Abel-Jacobi map has dimension
$ \ge m $, it is enough to take the
$ (m+1) $-ple fiber product of
$ X_{R} $ over
$ R $, and consider the pull-back of the above cycle by the projection to
the fiber product of the
$ i $-th and the
$ (m+1) $-th factors for
$ 1 \le i \le m $.

In the case
$ X $ is not necessarily projective, let
$ \rho : X' \rightarrow X $ be a birational morphism such that
$ X' $ is smooth projective.
We may assume that
$ X' $ is defined over
$ R $, and let
$ X'_{R}, X'_{K} $ be as above.
Then we have an isomorphism
$ \rho _{*} : L^{2}\CH_{0}(X')_{\Bbb Q} \rightarrow
L^{2}\CH_{0}(X)_{\Bbb Q} $ (see [19]) and a decomposition
$ H^{2n-2}(X',{\Bbb Q}_{\langle k \rangle}) = H^{2n-2}(X,{\Bbb Q}_{\langle
k \rangle})
\oplus N $, where
$ N $ has level zero.
So
$ N $ does not contribute to the
$ m $-dimensional vector space constructed in (4.2), and the assertion is
reduced to the projective case.

\medskip
\noindent
{\it Remarks.} (i) Let
$ X $ be a smooth hypersurface of
$ {\Bbb P}^{n+1} $ defined over
$ k $.
Then we have
$ R = K = k $ in the notation of (4.1).
Take a closed point
$ \Spec k' $ of
$ X_{k} $, and let
$$
\zeta = \Delta - [k':k]^{-1}(X_{k}\times _{k}\Spec k')
\in \CH^{n}(X_{k}\times _{k}X_{k})_{\Bbb Q},
$$
where
$ \Delta $ is the diagonal.
Let
$ K' $ be the function field of
$ X_{k} $, and
$ \zeta ' $ the restriction of
$ \zeta $ to
$ X_{K'} $ as before.
Assume that the mixed sheaves in (1.1) are mixed Hodge Modules.
Then
$$
\zeta ' \in L^{n}\CH^{n}(X_{K'})_{\Bbb Q}.
$$

Indeed, it is enough to show by (1.5)
$$
\Ext^{1}({\Bbb Q}_{k},H^{2i-1}(U/k,{\Bbb Q} (i)) =0,\quad \Hom({\Bbb
Q}_{k},H^{2j}(U/k,{\Bbb Q}(j)) = 0
$$
for
$ 0 < i \le n/2, 0 < j < n/2, $ and
$ U $ a sufficiently small open subvariety of
$ X_{k}, $ because
$ H^{r}(X/k,{\Bbb Q}) = {\Bbb Q}(-r/2) $ or zero for
$ r > n $.
But the assertion is clear by
$$
W_{2i-1}H^{2i-1}(U/k,{\Bbb Q}) = W_{2j}H^{2j}(U/k,{\Bbb Q}) = 0.
$$

We can show furthermore that
$ {\Gr}_{L}^{n}cl_{K'}(\zeta ') \ne 0 $, and the image of
$ {\Gr}_{L}^{n}cl $ is infinite dimensional as in (4.2) if
$ \deg X \ge n + 2 $.

(ii) Let
$ X, Y $ be smooth proper complex varieties of dimension
$ n, m $ respectively, and
$ \zeta \in C^{i}(Y,X) $ for a positive integer
$ i $ such that
$ i + m \le n $.
Assume they are defined over
$ k $ so that they come from
$ X_{k} $, etc.
Let
$ K $ be the function field of
$ Y_{k}/k $, and
$ \zeta_{K} $ the restriction of
$ \zeta_{k} $ to
$ X_{K} := X_{k}\otimes_{k}K $.
Assume further
$ \zeta_{K} \in L^{m}\CH^{i+m}(X_{K}) $.
Then
$ {\Gr}_{L}^{m}cl(\zeta_{K}) $ is nonzero if the image of
$ \zeta_{*} : H^{m}(Y, {\Bbb Q}) \rightarrow H^{2i+m}(X,{\Bbb Q})(i) $ has
level of Hodge structure
$ m $.

\proclaim{{\bf 4.5.~Lemma}}
Let
$ X, Y, Z $ be smooth proper
$ k $-varieties, and
$ \zeta \in C^{i}(X,Y), \zeta ' \in C^{j}(Y,Z) $.
Let
$ \zeta '' = \zeta \scirc \zeta ' \in C^{i+j}(X,Z) $.
Assume
$ X $ is irreducible, and let
$ K $ be the function field of
$ X $.
Put
$ Y_{K} = Y\otimes _{k}K $, and similarly for
$ Z_{K} $.
Let
$ \zeta _{K}, \zeta ''_{K} $ be the restrictions of
$ \zeta , \zeta '' $ to
$ Y_{K}, Z_{K} $, and
$ \zeta '_{K} \in C^{j}(Y_{K},Z_{K}) $ be the pull-back of
$ \zeta ' $.
Then
$$
\zeta ''_{K} = (\zeta '_{K})_{*}\zeta _{K}.
$$
\endproclaim

\demo\nofrills {Proof.\usualspace}
This follows from
$ \zeta \scirc \zeta ' = p_{*}i^{*}(\zeta \times \zeta ') $ where
$ i = id\times \delta \times id : X\times Y\times Z \rightarrow X\times
Y\times Y\times Z $ with
$ \delta $ the diagonal, and
$ p : X\times Y\times Z \rightarrow X\times Z $ is the projection.
\enddemo

\medskip
\noindent
{\bf 4.6.}
{\it Proof of} (0.4).
With the notation of (4.1), we have
$$
\Ext^{2p}({\Bbb Q}_{S/k_{R}}, \pi _{*}{\Bbb Q}_{X_{R}/k_{R}}(p)) =
\Ext^{2p}({\Bbb Q}_{k_{R}}, (a_{X_{R}/k_{R}})_{*}{\Bbb
Q}_{X_{R}/k_{R}}(p)),
$$
and
$ {H}_{\Cal D}^{2p}(X,{\Bbb Q}_{\langle k \rangle}(p)) $ is equal to
$$
\hbox{\rm ind\,lim}_{R} \Ext^{2p}({\Bbb Q}_{k_{R}},
(a_{X_{R}/k_{R}})_{*}{\Bbb Q}
_{X_{R}/k_{R}}(p)).
$$
So the first assertion is reduced to the injectivity of the cycle map
$$
\CH^{p}(X_{R})_{\Bbb Q} \rightarrow \Ext^{2p}({\Bbb Q}_{k_{R}},
(a_{X_{R}/k_{R}})_{*}{\Bbb Q}_{X_{R}/k_{R}}(p)).
$$
See (3.4.2).
Then it follows from the assumption by using the injective morphisms in
(1.3.5).
The last assertion is reduced to the following.

\proclaim{{\bf 4.7.~Proposition}}
Let
$ X $ be a smooth proper variety over
$ k $ such that
$ k $ is algebraically closed in the function field of
$ X $.
Let
$ Y $ be a divisor on
$ X $, and
$ \rho : X' \rightarrow X $ be a birational
$ k $-morphism such that
$ X' $ is smooth projective and
$ Y' := \rho ^{-1}(Y) $ is a divisor with normal crossings whose
irreducible components
$ Y'_{i} $ are smooth.
Then the Abel-Jacobi map over
$ k $ for
$ U := X \setminus Y $ and
$ p = 2 $ is injective if the same injectivity holds for
$ X' $ and
$ p = 2 $.
\endproclaim

\demo\nofrills {Proof.\usualspace}
By the compatibility of the cycle map with the pull-back, the assertion is
easily reduced to that for
$ U' = \rho ^{-1}(U) $.
So we assume
$ X' = X $ and
$ Y'_{i} $ will be denoted by
$ Y_{i} $.

Let
$ M = H^{3}(U/k,{\Bbb Q})(2) $.
Then it has weights
$ \ge - 1 $, and
$ \Ext^{i}({\Bbb Q}_{k},M) = \Ext^{i}({\Bbb Q}_{k},W_{0}M) $ for
$ i = 0, 1 $, where
$ W $ is the weight filtration.
So we have an exact sequence
$$
\Hom({\Bbb Q}_{k}, {\Gr}_{0}^{W}M) \rightarrow \Ext^{1}({\Bbb Q}_{k},
{\Gr}_{-1}^{W}
M) \rightarrow \Ext^{1}({\Bbb Q}_{k}, M).
\leqno(4.7.1)
$$
By the weight spectral sequence [12] we have furthermore
$$
{\Gr}_{-1}^{W}M = \Coker(\oplus _{i} H^{1}(Y_{i}/k,{\Bbb Q}(1)) \rightarrow
H^{3}(X/k,{\Bbb Q}(2)))
$$
and
$ {\Gr}_{0}^{W}M $ is isomorphic to a quotient of
$$
\Ker(\oplus _{i} H^{2}(Y_{i}/k,{\Bbb Q}(1)) \rightarrow H^{4}(X /k,{\Bbb
Q}(2))).
$$

Let
$ i : Y \rightarrow X $ denote the inclusion morphism.
Then we have
$$
\aligned
{\Gr}_{-1}^{W}M
&= \Coker(i_{*} : H^{3}(a_{Y/k})_{*}i^{!}{\Bbb Q}_{X/k}(2)
\rightarrow H^{3}(X/k,{\Bbb Q}(2))),
\\
M/W_{-1}M
&= \Ker(i_{*} : H^{4}(a_{Y/k})_{*}i^{!}{\Bbb Q}_{X/k}(2)
\rightarrow H^{4}(X/k,{\Bbb Q}(2))),
\endaligned
\leqno(4.7.2)
$$
using the associated long exact sequence, and the first morphism of (4.7.1)
is also induced by
$ (a_{Y/k})_{*}i^{!}{\Bbb Q}_{X/k} \rightarrow (a_{X/k})_{*}{\Bbb Q}_{X/k}
$.
So the cycle map gives a morphism of an exact sequence
$$
(\oplus _{i} \CH^{1}(Y_{i})_{\Bbb Q})^{(0)} \rightarrow L^{1}\CH
^{2}(X)_{\Bbb Q} \rightarrow L^{1}\CH^{2}(U)_{\Bbb Q}
\leqno(4.7.3)
$$
to (4.7.1), where
$ (\oplus _{i} \CH^{1}(Y_{i})_{\Bbb Q})^{(0)} = \Ker((\oplus _{i}
\CH^{1}(Y_{i})_{\Bbb Q}) \rightarrow H^{4}(X_{\Bbb C},{\Bbb Q}
(p))) $.
Since the Hodge conjecture for codimension one cycles on
$ Y_{i} $ is true, the morphism of the left terms is surjective.
See Remark (i) after (4.1).
This implies also the surjectivity of
$ L^{1}\CH^{2}(X)_{\Bbb Q} \rightarrow L^{1}\CH^{2}(U)_{\Bbb Q} $.
Then the assertion is reduced to the next proposition, because
$$
\Ext^{1}({\Bbb Q}_{k}, {\Gr}_{-1}^{W}M) = \Coker(\oplus _{i}
J^{1}(Y_{i}/k)_{\Bbb Q}
\rightarrow J^{2}(X/k)_{\Bbb Q})
\leqno(4.7.4)
$$
where
$ J^{p}(X/k)_{\Bbb Q} $ is as in (1.3.4).
\enddemo

\proclaim{{\bf 4.8.~Proposition}}
Let
$ X, Y $ be smooth projective
$ k $-varieties of pure dimension
$ n $ and
$ n - 1 $ respectively.
Let
$ g : Y \rightarrow X $ be a morphism of
$ k $-varieties.
Then
$$
{\Gr}_{L}^{1}cl(g_{*}L^{1}\CH^{1}(Y)_{\Bbb Q}) =
{\Gr}_{L}^{1}cl(L^{1}\CH^{2}(X)_{\Bbb Q})
\cap \Im(g_{*} : J^{1}(Y/k)_{\Bbb Q} \rightarrow J^{2}(X/k)_{\Bbb Q})
$$
\endproclaim

\demo\nofrills {Proof.\usualspace}
Let
$ l $ denote the hyperplane section class of
$ X $.
We will denote also by
$ l $ the induced morphism
$$
l: H^{i}(X/k,{\Bbb Q}) \rightarrow H^{i+2}(X/k,{\Bbb Q})(1).
$$
(It is defined over
$ k $ because the condition for a general hyperplane section is
Zariski-open, and
$ k $ is an infinite field.)
By the Lefschetz decomposition, we have
$$
H^{3}(X/k,{\Bbb Q}) = H^{3}(X/k,{\Bbb Q})^{\prim}
\oplus lH^{1}(X/k,{\Bbb Q}(-1)),
$$
where
$ H^{3}(X/k,{\Bbb Q})^{\prim} = \Ker \,l^{n-2} \subset H^{3}(X/k,{\Bbb Q})
$.
This implies
$$
J^{2}(X/k)_{\Bbb Q} = J^{2}(X/k)_{\Bbb Q}^{\prim} \oplus
lJ^{1}(X/k)_{\Bbb Q},
$$
where
$ J^{2}(X/k)_{\Bbb Q}^{\prim} =
\Ext^{1}({\Bbb Q}_{k}, H^{3}(X/k,{\Bbb Q}(2))^{\prim}) $.

Let
$ e \in J^{1}(Y/k)_{\Bbb Q} $ such that
$ g_{*}e \in {\Gr}_{L}^{1}cl(L^{1}\CH^{2}(X)_{\Bbb Q}) $.
We have to show
$ g_{*}e \in {\Gr}_{L}^{1}cl(g_{*}L^{1}\CH^{1}(Y)_{\Bbb Q}) $.
We first reduce the assertion to the case
$ g_{*}e \in J^{2}(X/k)_{\Bbb Q}^{\prim} $.

Since
$ l $ is algebraically defined as a correspondence, we have
$$
l^{n-2}g_{*}e \in {\Gr}_{L}^{1}cl(L^{1}\CH^{n}(X)_{\Bbb Q}) \subset
J^{n}(X/k)_{\Bbb Q}.
$$
By (3.6), we can identify
$ {\Gr}_{L}^{1}cl(L^{1}\CH^{n}(X)_{\Bbb Q}) $ with
$ \Alb_{X/k}(k)_{\Bbb Q} $.
Let
$ P_{Y/k} $ be the Picard variety of
$ Y $.
Since
$ l^{n-2}g_{*} $ is algebraically defined by a correspondence
$ \Gamma $, we have a commutative diagram
$$
\CD
P_{Y/k}(k)_{\Bbb Q}
@>>>
J^{1}(Y/k)_{\Bbb Q}
\\
@VV{\Gamma _{*}}V @VV{\Gamma _{*}}V
\\
\Alb_{X/k}(k)_{\Bbb Q} @>>> J^{n}(X/k)_{\Bbb Q}
\endCD
$$
where the horizontal morphisms are injective.
We see that the left vertical morphism comes from a morphism of abelian
varieties.
(Indeed, if
$ D $ denotes the Poincar\'e divisor on
$ P_{Y/k}\times _{k}Y $, the composition
$ D \scirc \Gamma \in \CH^{n}(P_{Y/k}\times _{k}X) $ as correspondences
gives further
$ (D \scirc \Gamma )_{\Alb} \in \CH^{m}(P_{Y/k}\times _{k}\Alb _{X/k}) $
which is given by the graph of a morphism of abelian varieties,
and induces
$ \Gamma _{*} : P_{Y/k}(k) \rightarrow \Alb_{X/k}(k) $, where
$ m = \dim \Alb_{X/k} $.)
Then the diagram induces injective morphisms on the kernel and the cokernel
of
$ \Gamma _{*} $ (by using the base change by
$ k \rightarrow {\Bbb C} $ together with (1.1.2)).
So
$ l^{n-2}g_{*}e $ belongs to
$ \Gamma _{*}P_{Y/k}(k)_{\Bbb Q} $, and there exists
$ \zeta \in L^{1}\CH^{1}(Y)_{\Bbb Q} $ such that
$$
l^{n-2}g_{*}e = l^{n-2}g_{*}{\Gr}_{L}^{1}cl(\zeta ).
$$
Then we may assume
$ l^{n-2}g_{*}e = 0 $ (i.e.
$ g_{*}e $ is primitive) by replacing
$ e $ with
$ e - {\Gr}_{L}^{1}cl(\zeta ) $.

Now let
$$
P_{Y/k}(k)_{\Bbb Q}^{(0)} = \Ker \,\Gamma _{*}
\subset P_{Y/k}(k)_{\Bbb Q},\,\,\,J^{1}(Y/k)_{\Bbb Q}^{(0)}
=\Ker \,\Gamma _{*} \subset J^{1}(Y/k)_{\Bbb Q},
$$
and consider
$$
g^{*}l^{n-3}g_{*}e \in \Alb_{Y/k}(k)_{\Bbb Q} \subset J^{n-1}(Y/k)_{{\Bbb
Q}}.
$$
By a similar argument, there exists
$ \zeta ' \in P_{Y/k}(k)_{\Bbb Q}^{(0)} $ such that
$$
g^{*}l^{n-3}g_{*}e = g^{*}l^{n-3}g_{*}{\Gr}_{L}^{1}cl(\zeta ').
$$
So it remains to show that the restriction of
$ g^{*}l^{n-3} $ to
$ \Im \, g_{*} \cap H^{3}(X/k,{\Bbb Q})^{\prim} $ is injective, or
equivalently, the pairing
$ \langle u,l^{n-3}v \rangle $ is nondegenerate on
$ \Im \, g_{*} \cap H^{3}(X/k,{\Bbb Q})^{\prim} $.
But this is clear because the pairing is a polarization of Hodge structure.
This completes the proof of (4.7-8) and (0.4).
\enddemo

\medskip
\noindent
{\it Remarks.} (i) We can verify that Bloch's conjecture [7] is true if the
conclusion of (0.4) holds.
Let
$ X $ be a smooth proper complex variety of pure dimension
$ 2 $, and
$ N^{1}H^{2}(X,{\Bbb Q}_{\langle k \rangle}) $ the space of Hodge (i.e.
algebraic) cycles of codimension one, which is a direct sum of
$ {\Bbb Q}_{\langle k \rangle}(-1) $.
Then for
$ \zeta \in L^{2}\CH^{2}(X)_{\Bbb Q} $, we can show that
$ cl(\zeta ) = 0 $ if
$ cl(\zeta ) \in \Ext^{2}({\Bbb Q}_{\langle k \rangle}, N^{1}H^{2}(X,{\Bbb
Q}
_{\langle k \rangle}(2))) $.

Indeed, let
$ C $ be a (not necessarily connected) smooth proper curve with a morphism
$ g : C \rightarrow X $ such that
$ N^{1}H^{2}(X,{\Bbb Q}_{\langle k \rangle}) = g_{*}H^{0}(C,{\Bbb
Q}_{\langle k \rangle}( -1)) $.
Then we have an injection
$$
g^{*} : N^{1}H^{2}(X,{\Bbb Q}_{\langle k \rangle}) \rightarrow
H^{2}(C,{\Bbb Q}
_{\langle k \rangle}),
$$
but
$ g^{*}cl(\zeta ) = cl(g^{*}\zeta ) = 0 $.
So
$ cl(\zeta ) = 0 $.
See also [40,~II], (4.12).

In general, it is easy to show that for a smooth projective variety
$ X $ and
$ \zeta \in L^{r}\CH^{p}(X)_{\Bbb Q} $ with
$ r > p $ we have
$ {\Gr}_{L}^{r}cl(\zeta ) = 0 $ and hence
$ cl(\zeta ) = 0 $ (using the pull-back to an intersection of generic
hyperplanes of dimension
$ 2p - r) $.
See [40,~I].

(ii) The extension groups in
$ \MHS({\Bbb Q})_{\langle k \rangle} $ are generally too big.
For example, we have
$$
\Ext^{2}({\Bbb Q}_{\langle k \rangle}, {\Bbb Q}_{\langle k \rangle}(1)) \ne
0
$$
by (1.4--5) in the case the mixed sheaves in (1.1) are mixed Hodge Modules.
Indeed, for
$ R, S $ as in (2.3), we have
$$
\Ext^{1}({\Bbb Q}_{k}, H^{1}(S/k_{R}, {\Bbb Q}(1))) \ne 0,
$$
if
$ {\Gr}_{1}^{W}H^{1}(S_{\Bbb C},{\Bbb Q}) \ne 0 $.

(iii) With the notation of (0.4), we have a commutative diagram
$$
\CD
L^{1}\CH^{p}(Y)_{\Bbb Q} @>>> J^{p}(Y/k)_{\Bbb Q}
\\
@VVV @VVV
\\
L^{1}\CH^{p}(Y_{\Bbb C})_{\Bbb Q} @>>>
J^{p}(Y_{\Bbb C})_{\Bbb Q}
\endCD
$$
where the horizontal morphisms are the Abel-Jacobi maps.
If the composition of the horizontal and vertical morphisms is injective,
it would imply the injectivity of the Abel-Jacobi map defined over
$ k $.
But the converse is not clear, because the right vertical morphism is not
injective by (1.4) (although so is the left one).

\medskip
\noindent
{\bf 4.9.}
{\it Remarks.} (i) Assume that a smooth proper complex variety
$ X $ admits the K\"unneth decomposition in the Chow group in the sense of
Murre [32] [33] (i.e., if there exist
$ p_{i} \in C^{0}(X,X) $ for
$ i \in {\Bbb Z} $ such that
$ p_{i} = 0 $ for
$ i > 2 \dim X $ or
$ i < 0 $,
$ \sum _{i} p_{i} $ is the diagonal of
$ X $,
$ p_{i}p_{j} = \delta _{i,j}p_{i} $, and the action of
$ p_{i} $ on
$ H^{j}(X,{\Bbb Q}) $ is
$ \delta _{i,j}id) $.
Let
$ q_{r,j} = \sum _{0\le i<j} p_{r - i} $ and
$ q'_{r,j} = \sum _{i<0} p_{r - i} + \sum _{i\ge j} p_{r -i} $.
Then, following Murre, we define the filtration
$ F_{M} $ on
$ \CH^{p}(X)_{\Bbb Q} $ by
$$
{F}_{M}^{j}\CH^{p}(X)_{\Bbb Q} = (q'_{2p,j})_{*}\CH^{p}(X)_{\Bbb Q} =
\Ker((q_{2p,j})_{*} : \CH^{p}(X)_{\Bbb Q}
\rightarrow \CH^{p}(X)_{\Bbb Q}).
$$

Let
$ L $ be as in (3.4.3).
Then we can verify
$$
F_{M} \subset L
\leqno(4.9.1)
$$
inductively, using
$ p_{i}p_{j} = \delta _{i,j}p_{i} $ together with (3.5).
Similarly, we see
$$
{\Gr}_{{F}_{M}}^{j}{\Gr}_{L}^{i} = 0\quad \text{for }i \ne j,
\leqno(4.9.2)
$$
because the action of
$ p_{2p - j} $ on
$ {\Gr}_{{F}_{M}}^{j} $ is the identity.
In particular,
$$
F_{M} = L \mod \bigcap _{i} L^{i}.
\leqno(4.9.3)
$$
If the cycle map (3.4.1) is injective for
$ X, $ we have
$ \bigcap _{i} L^{i} = 0 $ and the action of
$ p_{i} $ on
$ \CH^{p}(X)_{\Bbb Q} $ for
$ i > 2p $ or
$ i < p $ is zero.
(The relation with a conjectural filtration of Beilinson and Bloch is
explained in [26] assuming the algebraicity of the K\"unneth components of
the diagonal.)

(ii) In [43], Shuji Saito has defined a filtration
$ F_{\Sh} $ on
$ \CH^{p}(X)_{\Bbb Q} $ for a smooth projective complex variety
$ X $ by letting
$ {F}_{\Sh}^{0}\CH^{p}(X)_{\Bbb Q} = \CH^{p}(X)_{\Bbb Q} $ and using an
inductive formula
$$
{F}_{\Sh}^{j+1}\CH^{p}(X)_{\Bbb Q} = \sum \Gamma
_{*}{F}_{\Sh}^{j}\CH^{p-i}(Y)_{\Bbb Q}.
$$
Here the summation is taken over correspondences
$ \Gamma \in C^{i}(Y,X)_{\Bbb Q} $ with
$ i $ integers and
$ Y $ smooth projective complex varieties such that the image of
$$
\Gamma _{*} : H^{2p-j-2i}(Y,{\Bbb Q})(-i) \rightarrow H^{2p-j}(X,{\Bbb Q})
\leqno(4.9.4)
$$
is contained in
$ N^{p-j+1}H^{2p-j}(X,{\Bbb Q}) $.
(Here
$ N $ denotes the ``coniveau'' filtration.)
It is possible to consider some variants of this definition.
We get filtrations
$ F_{\alpha }, F_{\beta }, F_{\gamma } $ by modifying the condition on
$ \Gamma $ respectively as follows.

\noindent
($\alpha $) There exists
$ \Gamma ' \in C^{i}(Y,X)_{\Bbb Q} $ such that the morphisms (4.9.4) for
$ \Gamma $ and
$ \Gamma ' $ coincide and
$ \Gamma ' $ is supported on
$ Y \times Z $ with
$ \dim Z < \dim X - p + j $.

\noindent
($\beta $) The morphism (4.9.4) is zero.

\noindent
($\gamma $) The morphism (4.9.4) is zero and
$ Y = X, i = 0 $.

We have clearly
$ F_{\Sh} \supset F_{\alpha } \supset F_{\beta } \supset F _{\gamma } $
(and they contain
$ F_{M} $ if the latter exists).
The filtrations
$ F_{\alpha }, F_{\beta } $ are stable by the action of correspondences
(using [19] for
$ F_{\alpha }) $.
We can show also
$$
F_{\alpha } \subset L.
\leqno(4.9.5)
$$
But it does not seem easy to verify
$ F_{\Sh} = F_{\alpha } $ or
$ F_{\Sh} \subset L $ without assuming the standard conjectures [28].
If the K\"unneth components of the diagonal are algebraic for
$ X $, then
$$
{\Gr}_{{F}_{*}}^{j}{\Gr}_{L}^{i} = 0\quad \text{for }i \ne j\,\,\,
\text{and }* = \Sh, \alpha , \beta , \gamma.
\leqno(4.9.6)
$$
So we have
$ F_{*} = L $ for
$ * = \alpha , \beta , \gamma $, if the K\"unneth components of the
diagonal are algebraic and the cycle map (3.4.1) is injective for
$ X $.

\bigskip\bigskip
\centerline{{\bf 5.~Higher Chow Groups}}

\bigskip
\noindent
{\bf 5.1.}
Let
$ k $ be a subfield of
$ {\Bbb C} $ as above.
For a
$ k $-variety
$ X $ of pure dimension
$ n $, let
$ \CH^{p}(X,m)_{\Bbb Q} $ denote Bloch's higher Chow group with rational
coefficients [8].
Its elements are represented by cycles of codimension
$ p $ on
$ X \times _{k} {\Bbb A}^{m} $.
Let
$ {\Cal M}(X/k,{\Bbb Q}) $ be the categories of mixed sheaves for
$ k $-varieties
$ X $ as in (1.1).
Let
$ d = n - p $.
Then by [37] we can construct a cycle map
$$
cl : \CH^{p}(X,m)_{\Bbb Q} \rightarrow {H}_{2d+m}^{\Cal D}(X/k,{\Bbb
Q}(d)),
\leqno(5.1.1)
$$
where the target is defined to be
$$
\Ext^{-2d-m}({a}_{X/k}^{*}{\Bbb Q}_{k}, {a}_{X/k}^{!}{\Bbb Q}_{k}(-d))
(= \Ext^{-2d-m}({\Bbb Q}_{k}, (a_{X/k})_{*}{a}_{X/k}^{!}{\Bbb Q}_{k}(-d))),
$$
with
$ \Ext $ taken in
$ D^{b}{\Cal M}(X/k,{\Bbb Q}) $ (or
$ D^{b}{\Cal M}(\Spec k/k,{\Bbb Q})) $.
If
$ X $ is smooth, it becomes
$$
\aligned
{H}_{\Cal D}^{2p-m}(X/k,{\Bbb Q}(p)):= \Ext^{2p-m}
&({\Bbb Q}_{X/k}, {\Bbb Q}_{X/k}(p))
\\
(= \Ext^{2p-m}
&({\Bbb Q}_{k}, (a_{X/k})_{*}{\Bbb Q}_{X/k}(p))),
\endaligned
$$
in the notation of (1.3).
If furthermore
$ X $ is smooth proper and
$ m \ge 1 $, then the cycle map induces a generalized Abel-Jacobi map
$$
cl' : \CH^{p}(X,m)_{\Bbb Q} \rightarrow \Ext^{1}({\Bbb Q}_{k},
H^{2p-m-1}(X/k, {\Bbb Q}(p))),
\leqno(5.1.2)
$$
by taking the composition with the canonical morphism
$$
{H}_{\Cal D}^{2p-m}(X/k,{\Bbb Q}(p)) \rightarrow \Ext^{1}({\Bbb Q}_{k},
H^{2p-m-1}(X /k, {\Bbb Q}(p))),
$$
because
$ H^{2p-m}(X/k, {\Bbb Q}(p)) $ is pure of weight
$ - m $.

Let now
$ X $ be a complex algebraic variety.
Then (5.1.1) induces the cycle map
$$
cl : \CH^{p}(X,m)_{\Bbb Q} \rightarrow {H}_{2d+m}^{\Cal D}(X,{\Bbb
Q}_{\langle k \rangle}(d)),
\leqno(5.1.3)
$$
where the target is defined in a similar way by replacing
$ {\Cal M}(X/k,{\Bbb Q}) $ with
$ {\Cal M}(X,{\Bbb Q})_{\langle k \rangle} $.
If
$ X $ is smooth, it becomes
$ {H}_{\Cal D}^{2p-m}(X,{\Bbb Q}_{\langle k \rangle}(p)) $ in (3.3).

If
$ m = 1 $, we have a natural map
$$
\CH^{p-1}(X)_{\Bbb Q}\otimes _{{\Bbb Z}}{\Bbb C}^{*} \rightarrow
\CH^{p}(X,1)_{\Bbb Q},
\leqno(5.1.4)
$$
using the isomorphism
$ \CH^{1}(pt,1) = {\Bbb C}^{*} $.
See [8].
(Indeed, an element of
$ \CH^{p}(X,1) $ is represented by
$ \sum _{i} (Z_{i},g_{i}) $ where
$ Z_{i} $ is an integral closed subvariety of codimension
$ p -1 $ in
$ X $ and
$ g_{i} $ are nonzero rational functions on
$ Z_{i} $ such that
$ \sum _{i} \div g_{i} = 0 $ in
$ X $.)
Note that (5.1.4) is not injective.
The image of (5.1.4) is denoted by
$ \CH_{\dec}^{p}(X,1)_{\Bbb Q} $, and its elements are called
decomposable.
We denote the cokernel of (5.1.4) by
$ \CH_{\ind}^{p}(X,1)_{\Bbb Q} $.
See [20], [30].

Consider the cycle map to the usual Deligne cohomology for a smooth
projective variety
$ X $ over
$ {\Bbb C} $
$$
cl : \CH^{p}(X,1)_{\Bbb Q} \rightarrow {H}_{\Cal D}^{2p-1}(X,{\Bbb
Q}(p)) = {\Ext}_{\MHS}^{1}
({\Bbb Q}, H^{2p-2}(X, {\Bbb Q}(p))).
\leqno(5.1.5)
$$
See [2], [9], [16], etc.
The image of
$ \CH_{\dec}^{p}(X,1)_{\Bbb Q} $ by (5.1.5) coincides with
$$
{\Ext}_{\MHS}^{1}({\Bbb Q}, N^{p-1}H^{2p-2}(X, {\Bbb Q}(p))) =
N^{p-1}H^{2p-2}(X, {\Bbb Q}(p-1))\otimes _{{\Bbb Z}}{\Bbb C}^{*},
$$
where
$ N^{p-1}H^{2p-2}(X, {\Bbb Q}(p-1)) $ is the subspace of algebraic cycle
classes (which is contained in the subspace of Hodge cycles).
So we cannot detect the image of
$ {\CH}_{\hom}^{p-1}(X)_{\Bbb Q}\otimes _{{\Bbb Z}}{\Bbb C}^{*} $ in
$ \CH^{p}(X,1)_{\Bbb Q} $ by using (5.1.5),
$ $ where
$ {\CH}_{\hom}^{p-1}(X)_{\Bbb Q} $ denotes the subgroup of homologically
equivalent to zero cycles.

\proclaim{{\bf 5.2.~Theorem}}
Let
$ X $ be a smooth complex projective variety such that
$ N^{p-2}H^{2p-3}(X,{\Bbb Q}) \ne 0 $, where
$ N $ is the coniveau filtration.
Then for any positive integer
$ m $, there exist
$ \zeta _{i} \in {\CH}_{\hom}^{p-1}(X)_{\Bbb Q}
\,(1 \le i \le m) $ and a subfield
$ K $ of
$ {\Bbb C} $ finitely generated over
$ k $ such that for any complex numbers
$ \alpha _{1}, \dots , \alpha _{m} $ not algebraic over
$ K $, the images of
$ \zeta _{i}\otimes \alpha _{i}\,(1 \le i \le m) $ by the composition of
(5.1.4) and (5.1.3) are linearly independent over
$ {\Bbb Q} $.
\endproclaim

\noindent
{\it Remark.} If
$ p = 2 $ or
$ \dim X + 1 $, the assumption on the coniveau filtration is equivalent to
the nonvanishing of
$ H^{1}(X,{\Cal O}_{X}) $.
Furthermore, (5.1) in the case
$ p = \dim X + 1 $ holds by replacing that condition with the nonvanishing
of
$ H^{2}(X,{\Cal O}_{X}) $ using an argument similar to (4.4).

\medskip\noindent
{\it Proof of} (5.2).
By assumption on the coniveau filtration, there exist an irreducible smooth
complex projective curve
$ C $ and
$ \zeta ' \in \CH^{p-1}(C\times X) $ such that the induced morphism
$ H^{i}(C,{\Bbb Q}) \rightarrow H^{2p-4+i}(X,{\Bbb Q}(p-2)) $ is
nontrivial for
$ i = 1 $ and vanishes for
$ i = 2 $.
(Indeed, for a closed subvariety
$ Y $ of codimension
$ p - 2 $ in
$ X $, there is a curve together with a divisor on
$ C\times \widetilde{Y} $ such that the induced morphism
$ H^{i}(C,{\Bbb Q}) \rightarrow H^{i}(\widetilde{Y},{\Bbb Q}) $ is surjective
for
$ i = 1 $ and vanishes for
$ i = 2 $, where
$ \widetilde{Y} \rightarrow Y $ is a resolution of singularities.)

Let
$ R $,
$ S $,
$ X_{R} $ be as in (2.1).
We may assume that
$ C $ and
$ \zeta ' $ are defined over
$ R $, i.e. there exist an
$ R $-scheme
$ C_{R} $ and
$ \zeta _{R} \in \CH^{p-1}(C_{R}\times _{R}X_{R}) $ such that their base
changes by
$ R \rightarrow {\Bbb C} $ are
$ C $ and
$ \zeta ' $.
Let
$ U $ be a nonempty affine open subvariety of
$ C_{R} $.
We define
$ S' $ to be the fiber product of
$ m $ copies of
$ U $ over
$ S = \Spec R $, and
$ \zeta '_{i,S'} \in \CH^{p-1}(S'\times_{R}X_{R})$ to be the pull-back of
$ \zeta '_{R} $ by
$ p_{i}\times id $ where
$ p_{i} $ is the composition of the projection to the
$ m $-th factor
$ S' \rightarrow U $ and the inclusion
$ U \rightarrow C_{R} $.
Let
$ R' $ be the affine ring of
$ S' $.
Then
$ R' $ is integral.
We choose an embedding of
$ K = k(S') $ into
$ {\Bbb C} $.
The condition on
$ \alpha _{i} $ is relevant to this
$ K $ (i.e. the
$ \alpha _{i} $ are not algebraic over
$ K) $.
Let
$ \zeta_{i} \in \CH^{p-1}(X) $ be the base change of
$ \zeta _{i,S'} $.
It is homologically equivalent to zero by the assumption on
$ \zeta ' $.

Let
$ R'' $ be a finitely generated smooth
$ k $-subalgebra of
$ {\Bbb C} $ such that
$ R' \subset R'' $,
$ \Spec R'' $ is smooth over
$ \Spec R' $, and
$ \alpha _{i} \in R'' $ for any
$ i $.
Let
$ k_{R} $ be the algebraic closure of
$ k $ in
$ R $.
We may assume that
$ k_{R} $ is also algebraically closed in
$ R'' $ by replacing
$ R $,
$ R' $ with the base change by
$ k_{R} \rightarrow k_{R''} $ (and
$ k_{R} $ with
$ k_{R''} $).
Note that the algebraic closure of
$ K $ is unchanged.
Let
$ \pi : S'' \rightarrow S' $ be a smooth projective compactification of
$ \Spec R'' \rightarrow S' $.
Let
$ D_{i} $ be the support of the divisor defined by
$ \alpha _{i} $.
Put
$ D = \cup _{i} D_{i} $.
Note that
$ D_{i} $ are nonempty by assumption on
$ \alpha _{i} $.

We have to show that the images of
$ \zeta_{i}\otimes \alpha_{i} $ by the cycle map are linearly independent.
Let
$ U'' $ be a nonempty open subvariety of
$ S'' $.
We may assume that
$ U''\cap D $ is nonempty and smooth over
$ S' $, because it is only required that any nonempty open subvariety of
$ \Spec R'' $ contains
$ U''\setminus D $ for some
$ U'' $.
We can prove a commutative diagram
$$
\CD
\CH^{p}((U''\setminus D)\times _{R}X_{R},1)_{\Bbb Q}
@>>>
\CH^{p-1}((U''\cap D)\times _{R}X_{R})_{\Bbb Q}
\\
@VVV @VVV
\\
{H}_{\Cal D}^{2p-1}((U''\setminus D)\times _{R}X_{R}/k_{R},{\Bbb Q}
(p))
@>>>
{H}_{\Cal D}^{2p-2}((U''\cap D)\times _{R}X_{R}/k_{R},{\Bbb Q}(p-1))
\endCD
$$
See for example [42].
By definition
$ \zeta _{i}\otimes \alpha _{i} $ is represented by the pull-back of
$ \zeta _{i,S'} $ to
$ (U''\setminus D)\times _{R}X_{R} $ together with the function
$ \alpha _{i} $ coming from
$ U''\setminus D $.
It is enough to show that the images of their nontrivial linear
combinations in
$ {H}_{\Cal D}^{2p-2}((U''\cap D)\times _{R}X_{R},{\Bbb Q}(p)) $ are
nonzero by the same argument as in (4.2) using (2.6) together with (4.1.2)
applied to
$ (U''\setminus D)\times _{R}X_{R}
\rightarrow U''\setminus D $.
Consider the image of the cycle in
$ \CH^{p-1}((U''\cap D)\times _{R}X_{R})_{\Bbb Q} $.
It is identified with a linear combination of the pull-back of
$ \zeta _{i,S'} $ by
$ (U''\cap D)\times _{R}X_{R} \rightarrow S'\times _{R}X_{R} $ where the
coefficients are multiplied by the multiplicity of zero or pole of
$ \alpha _{i} $ along the divisor.
So the assertion follows from the same argument as in (4.2).

\proclaim{{\bf 5.3.~Theorem}}
Let
$ X $ be a smooth projective complex algebraic variety.
Assume
$ k $ is a number field.
Then the cycle map (5.1.3) for
$ p = 2, m = 1 $ is injective if the generalized Abel-Jacobi map (5.1.2)
for the same
$ p, m $ is injective for any
$ k $-smooth projective models of
$ X $ (where the notion of
$ k $-smooth projective model is defined as in (0.4)).
\endproclaim

\demo\nofrills {Proof.\usualspace}
Let
$ \zeta' \in \CH^{2}(X,1)_{\Bbb Q} $.
Replacing
$ R $ if necessary, we may assume that
$ \zeta' $ is defined over
$ R $, i.e.
there exists a cycle
$ \zeta' _{R} $ on
$ X_{R} \times _{R} {\Bbb A}_{R}^{1} $ such that
$ \zeta' _{R} \otimes _{R} {\Bbb C} = \zeta' $.
Let
$ Y = X_{R} $.
Then there exist smooth projective compactifications
$ \overline{Y} $,
$ \overline{S} $ of
$ Y, S $ together with a projective morphism
$ \overline{Y} \rightarrow \overline{S} $ extending
$ X_{R} \rightarrow S $.
Let
$ Z = $
$ \overline{Y} \setminus Y $.
We may assume it is a divisor on
$ \overline{Y} $.
Let
$ k_{R} $ be the algebraic closure of
$ k $ in
$ k(S) $.
Then
$ k_{R} $ is also algebraically closed in
$ k(Y) $ (because the generic fiber of
$ Y \rightarrow S $ is geometrically irreducible.)
We may assume that each irreducible components of
$ Z $ are geometrically irreducible (taking the base change of
$ \overline{Y} $,
$ S $ by a finite extension of
$ k_{R} $ if necessary).

Let
$ n = \dim Y $.
Then the cycle map induces a morphism of the localization sequence in [8]
to the corresponding exact sequence (see e.g. [42]):
$$
\CD
\CH^{1}(Z,1)_{\Bbb Q} @>>> \CH^{2}(\overline{Y},1)_{\Bbb Q}
@>>> \CH^{2}(Y,1)_{\Bbb Q} @>>> \CH^{1}(Z)_{\Bbb Q}
\\
@VVV @VVV @VVV @VV{(*)}V
\\
{H}_{2n-3}^{\Cal D}(Z,n-2) @>>> {H}_{\Cal D}^{3}(\overline{Y},2)
@>>> {H}_{\Cal D}^{3}(Y,2) @>>> {H}_{2n-4}^{\Cal D}(Z,n-2)
\endCD
$$
where
$ {H}_{2n-3}^{\Cal D}(Z,n-2) $ is the abbreviation of
$ {H}_{2n-3}^{\Cal D}(Z/k_{R},{\Bbb Q}(n-2)), $ etc.

Let
$ \zeta = \zeta'_{R} \in \CH^{2}(Y,1)_{\Bbb Q} $, and assume
$ cl(\zeta ) = 0 $.
Since
$ (*) $ is injective (using the base change by
$ k_{R} \rightarrow {\Bbb C} $, see loc.~cit),
$ $ we see that
$ \zeta $ is the image of some
$ \overline{\zeta } \in \CH^{2}(\overline{Y},1)_{\Bbb Q} $.
Furthermore,
$ cl(\overline{\zeta }) $ comes from
$ {H}_{2n-3}^{\Cal D}(Z/k_{R},{\Bbb Q}(n-2)) $ by the commutativity.
Since (5.1.2) is injective by hypothesis, it is enough to show that the
image of
$ \overline{\zeta } $ in
$ \Ext^{1}({\Bbb Q}_{k_{R}},H^{2}(\overline{Y}/k_{R},{\Bbb Q}(2))) $ comes
from
$ \CH^{1}(Z,1)_{\Bbb Q} $.
So the assertion is reduced to the following:

\proclaim{{\bf 5.4.~Proposition}}
Let
$ Y $ be a smooth projective
$ k $-variety of pure dimension
$ n $, and
$ Z $ a divisor on it.
Assume
$ Y $ and the irreducible components of
$ Z $ are geometrically irreducible over
$ k $.
Then
$$
\gathered
cl'(\CH^{2}(Y,1)_{\Bbb Q}) \cap \Im({H}_{2n-3}^{\Cal D}
(Z/k,{\Bbb Q}(n-2)) \rightarrow
\Ext^{1}({\Bbb Q}_{k},H^{2}(Y/k,{\Bbb Q}(2))))
\\
= \Im(\CH^{1}(Z,1)_{\Bbb Q} \rightarrow \Ext^{1}({\Bbb Q}_{k},
H^{2}(Y/k,{\Bbb Q}(2)))).
\endgathered
$$
\endproclaim

\noindent
{\it Proof.} Let
$ \widetilde{Z} \rightarrow Z $ be a resolution of singularities.
By Lemma (5.5) below, we may replace
$ Z $ in the formula by
$ \widetilde{Z} $.
Then
$ {H}_{2n-3}^{\Cal D}(Z/k,{\Bbb Q}(n-2)) $ becomes
$ \Ext^{1}({\Bbb Q}_{k},H^{0}(\widetilde{Z}/k,{\Bbb Q}(1))) $, and the
morphism of Deligne cohomologies is induced by the Gysin morphism for
$ g : \widetilde{Z} \rightarrow Y $.
Let
$ \xi = cl'(\zeta ) $ be an element of the left-hand side (with
$ Z $ replaced by
$ \widetilde{Z}) $.
We have the decomposition
$ \xi = \xi _{1} + \xi _{2} $ corresponding to the primitive decomposition
as in (4.8), where
$ \xi _{1} $ corresponds to the primitive part.

We first reduced the assertion to the case
$ \xi _{2} = 0 $.
Let
$ Z_{i} $ be the irreducible components of
$ Z $.
We may assume that the cohomology class of some
$ Z_{i} $ is not contained in the primitive part (because
$ \xi _{2} = 0 $ otherwise).
Let
$ l $ denote the hyperplane section class as in (4.8).
Consider the pushforward of
$ l^{n-1}\zeta $
$ \in \CH^{n+1}(Y,1)_{\Bbb Q} $ by
$ Y \rightarrow \Spec k $.
Let
$ \zeta _{i} \in \CH^{1}(\widetilde{Z}_{i},1)_{\Bbb Q} $ be its pull-back
by
$ \widetilde{Z}_{i} \rightarrow \Spec k $.
Then, modifying
$ \xi $ by a constant multiple of the image of
$ \zeta _{i} $, we may assume
$ \xi _{2} = 0 $.

Now let
$ V $ be a maximal dimensional vector subspace of
$ H^{0}(\widetilde{Z},{\Bbb Q}) $ which is sent injectively to a
subspace of the primitive part of
$ H^{2}(Y,{\Bbb Q}(1)) $ by the Gysin morphism for
$ g : \widetilde{Z} \rightarrow Y $.
Then there exists
$ \Gamma \in \CH^{0}(\widetilde{Z}\times _{k}\widetilde{Z})_{\Bbb Q} $ such
that
the action of the composition
$ \Gamma _{*}\scirc g^{*}\scirc l^{n-2}\scirc g_{*} $ on
$ V $ is the identity.
So
$ \xi $ comes from
$ (\Gamma _{*}\scirc g^{*}\scirc l^{n-2})\zeta \in
\CH^{1}(\widetilde{Z},1)_{\Bbb Q} $, and the assertion follows from the
compatibility of the cycle map with
$ \Gamma _{*} $,
$ g^{*} $,
$ l^{n-2} $,
$ g_{*} $.

\proclaim{{\bf 5.5.~Lemma}}
Let
$ X $ be a purely
$ n $-dimensional
$ k $-variety, and
$ \pi : \widetilde{X} \rightarrow X $ be a resolution of singularities.
Then
$$
{H}_{2n-1}^{\Cal D}(X/k,{\Bbb Q}(n-1)) = cl(\CH^{1}(X,1)_{\Bbb Q}) +
\pi _{*}({H}_{2n-1}^{\Cal D}(\widetilde{X}/k,{\Bbb Q}(n-1)))
$$
\endproclaim

\noindent
{\it Proof.} Let
$ U $ be a dense open subvariety of
$ X $ such that
$ \widetilde{U} := \pi ^{-1}(U) \rightarrow U $ is an isomorphism and
$ Z := X \setminus U $ and
$ \widetilde{Z} := \widetilde{X} \setminus \widetilde{U} $ are divisors.
Let
$ \xi \in {H}_{2n-1}^{\Cal D}(X/k,{\Bbb Q}(n-1)) $, and
$ \xi ' $ be its restriction to
$ \widetilde{U} = U $ which belongs to
$ {H}_{\Cal D}^{1}(U/k,{\Bbb Q}(1)) $.
Consider the commutative diagram induced by the cycle map:
$$
\CD
\CH^{1}(\widetilde{X},1)_{\Bbb Q} @>>> \CH^{1}(\widetilde{U},1)_{\Bbb Q}
@>>> \CH^{0}(\widetilde{Z})_{\Bbb Q} @>>> \CH^{1}(\widetilde{X})_{\Bbb Q}
\\
@VVV @VVV @VV{(*)}V @VV{(**)}V
\\
{H}_{2n-1}^{\Cal D}(\widetilde{X},n-1) @>>> {H}_{\Cal D}^{1}(\widetilde{U},1)
@>>> {H}_{2n-2}^{\Cal D}(\widetilde{Z},n-1)
@>>> {H}_{2n-2}^{\Cal D}(\widetilde{X},n-1)
\endCD
$$
where
$ {H}_{2n-1}^{\Cal D}(\widetilde{X},n-1) $ means
$ {H}_{2n-1}^{\Cal D}(\widetilde{X}/k,{\Bbb Q}(n-1)) $ as before.
By Lemma below,
$ (*) $ is an isomorphism and
$ (**) $ is injective (using the base change by
$ k \rightarrow {\Bbb C}) $.
So there exists
$ \zeta \in \CH^{1}(\widetilde{U},1)_{\Bbb Q} $ such that
$ \xi ' - cl(\zeta ) $ comes from an element
$ \xi '' $ of
$ {H}_{2n-1}^{\Cal D}(\widetilde{X}/k,{\Bbb Q}(n-1)) $.

Using a similar diagram with
$ \widetilde{X} $,
$ \widetilde{U} $,
$ \widetilde{Z} $ replaced by
$ X, U, Z $ and also the morphism between the two diagrams induced by the
direct image, we see that the image of
$ \zeta $ in
$ \CH^{0}(Z)_{\Bbb Q} $ is zero, because the image of
$ \zeta ' $ in
$ {H}_{2n-2}^{\Cal D}(Z/k,{\Bbb Q}(n-1)) $ is zero and
$ (*) $ remains injective with
$ \widetilde{Z} $ replaced by
$ Z $.
So
$ \zeta $ comes from an element of
$ \CH^{1}(X,1)_{\Bbb Q} $, and we get the assertion using the injectivity
of
$ {H}_{\Cal D}^{1}(X/k,{\Bbb Q}(1)) \rightarrow {H}_{\Cal D}^{1}(U/k,{\Bbb
Q}(1)) $, which follows from the vanishing of
$ {H}_{2n-1}^{\Cal D}(Z/k,{\Bbb Q}(n-1)) $.

\proclaim{{\bf 5.6.~Lemma}}
Let
$ X $ be a purely
$ n $-dimensional
$ k $-variety with irreducible components
$ X_{i} $.
Then we have a canonical isomorphism
$$
{H}_{2n}^{\Cal D}(X,{\Bbb Q}(n))
= \oplus _{i} {H}_{2n}^{\Cal D}(X_{i},{\Bbb Q}(n))
= \oplus _{i} {\Bbb Q}.
$$
\endproclaim

\noindent
{\it Proof.} By the localization sequence, we may assume
$ X $ smooth by deleting a closed subvariety of dimension
$ < n $ in
$ X $.
So the first isomorphism is clear and we may assume further
$ X $ irreducible.
Then the assertion becomes
$$
\Hom({\Bbb Q}_{X/k}, {\Bbb Q}_{X/k}) = {\Bbb Q},
$$
and it is clear if
$ X $ is geometrically irreducible.
In general, if the left-hand side is not
$ {\Bbb Q} $, there exists an endomorphism of
$ {\Bbb Q}_{X/k} $ inducing zero on some irreducible components of
$ X\otimes _{k}{\Bbb C} $ and nonzero on some other components (because the
endomorphisms on each component are the multiplications by rational
numbers).
But this is impossible because it comes from an endomorphism of
$ {\Cal O}_{X} $ as
$ {\Cal D}_{X} $-Module.
So the assertion follows.

\proclaim{{\bf 5.7.~Theorem}}
Let
$ X $ be a smooth complex projective variety.
If the cycle map (5.1.3) for
$ X $ is injective for
$ p = 2, m = 1 $, then it induces an injective morphism
$$
{\CH}_{\ind}^{2}(X,1)_{\Bbb Q} \rightarrow
\Ext^{1}({\Bbb Q}_{\langle k \rangle}, H^{2}(X,{\Bbb Q}_{\langle k
\rangle}(2))/ N^{1}H^{2}(X,{\Bbb Q}_{\langle k \rangle}(2))).
\leqno(5.7.1)
$$
\endproclaim

\noindent
{\it Remark.} It is not clear whether the morphism
$$
{\CH}_{\ind}^{2}(X,1)_{\Bbb Q} \rightarrow
\Ext_{\MHS}^{1}({\Bbb Q}, H^{2}(X,{\Bbb Q}(2))
/N^{1}H^{2}(X,{\Bbb Q}(2)))
\leqno(5.7.2)
$$
induced by (5.1.5) is injective.
Indeed, we have a natural morphism of the
target of (5.7.1) to that of (5.7.2), but the injectivity of its
restriction to the image of (5.7.1) seems to be nontrivial.
Note that (5.7) suggests that
$ {\CH}_{\ind}^{2}(X,1)_{\Bbb Q} = 0 $ if
$ H^{2}(X,{\Cal O}_{X}) = 0 $.
In the surface case, this is closely
related to Bloch's conjecture [7] due to [10].
See also [18].

\medskip\noindent
{\it Proof of} (5.7).
The assertion means that the inverse image of
$ \Ext^{1}({\Bbb Q}_{\langle k \rangle}, N^{1}H^{2}(X,{\Bbb Q}_{\langle k
\rangle}(2))) $ by the generalized Abel-Jacobi map
$$
cl' : \CH^{2}(X,1)_{\Bbb Q} \rightarrow
\Ext^{1}({\Bbb Q}_{\langle k \rangle}, H^{2}(X,{\Bbb Q}_{\langle k
\rangle}(2)))
$$
coincides with
$ {\CH}_{\dec}^{2}(X,1)_{\Bbb Q} $.
By definition
$ N^{1}H^{2}(X,{\Bbb Q}_{\langle k \rangle}) $ is isomorphic to a direct
sum of copies of
$ {\Bbb Q}_{\langle k \rangle}(-1) $ (generated by cycle classes), and is a
direct factor of
$ H^{2}(X,{\Bbb Q}_{\langle k \rangle}) $ by semisimplicity.
So we get
$$
cl'({\CH}_{\dec}^{2}(X,1)_{\Bbb Q}) = cl'({\CH}^{2}(X,1)_{\Bbb Q}) \cap
\Ext^{1}({\Bbb Q}_{\langle k \rangle}, N^{1}H^{2}(X,{\Bbb Q}_{\langle k
\rangle}(2))),
\leqno(5.7.3)
$$
by an argument similar to the proof of (5.4) using the compatibility of the
cycle map with correspondences.
(Note that the left-hand side is contained
in the right-hand side due to the compatibility of the cycle map with the
direct image by a morphism
$ Y \rightarrow X $, where
$ Y $ is a resolution of singularities of a divisor on
$ X $.)
Then the assertion is reduced to
$$
L^{2}\CH^{2}(X,1)_{\Bbb Q} \subset {\CH}_{\dec}^{2}(X,1)_{\Bbb Q},
\leqno(5.7.4)
$$
where the filtration
$ L $ on
$ \CH^{2}(X,1)_{\Bbb Q} $ is induced by
$ L $ on
$ {H}_{D}^{3}(X,{\Bbb Q}_{\langle k \rangle}(2)) $ in (3.3.3) using the
cycle map.
We have
$ {\Gr}_{L}^{3}\CH^{2}(X,1)_{\Bbb Q} = 0 $ due to the compatibility of
the cycle map with the pull-back by the inclusion of a point of
$ X $.
So the assertion is further reduced to the next proposition:

\proclaim{{\bf 5.8.~Proposition}}
With the above notation,
$ {\Gr}_{L}^{2}\CH^{2}(X,1)_{\Bbb Q} $ is generated by decomposable
higher cycles.
\endproclaim

\demo\nofrills {Proof.\usualspace}
Let
$ C $ be a smooth curve on
$ X $ obtained by intersecting general hyperplane sections of
$ X $.
Let
$ i : C \rightarrow X $ denote the inclusion.
Then there exists
$ \Gamma \in \CH^{1}(C\times X) $ such that the composition of
$ i^{*} $ and
$ \Gamma _{*} $ induces the identity on
$ H^{1}(X,{\Bbb Q}) $.
Since
$ \CH^{2}(C,1) $ consists of decomposable higher cycles, the assertion
follows from the compatibility of the cycle map with
$ i^{*} $ and
$ \Gamma _{*} $.
This completes the proofs of (5.8) and (5.7).
\enddemo

The following is a variant of a rigidity argument of Beilinson [2] and
M\"uller-Stach [30]:

\proclaim{{\bf 5.9.~Proposition}}
For a smooth complex projective variety
$ X $, let
$ \widetilde{N}^{i}H^{2i}(X,{\Bbb Q}_{\langle k \rangle}) $ denote the maximal
subobject of
$ H^{2i}(X,{\Bbb Q}_{\langle k \rangle}) $ which is isomorphic to a direct
sum of copies of
$ {\Bbb Q}_{\langle k \rangle}(-i) $.
Then the image of the morphism
$$
\CH^{p}(X,1)_{\Bbb Q} \rightarrow
\Ext^{1}({\Bbb Q}_{\langle k \rangle}, H^{2p-2}(X,{\Bbb Q}_{\langle k
\rangle}(p))
/\widetilde{N}^{p-1}H^{2p-2}(X,{\Bbb Q}_{\langle k \rangle}(p)))
\leqno(5.9.1)
$$
induced by the cycle map is countable.
\endproclaim

\demo\nofrills {Proof.\usualspace}
Let
$ {\Cal M} = H^{2p-2}(X,{\Bbb Q}_{\langle k \rangle}(p))
/\widetilde{N}^{p-1}H^{2p-2}(X,{\Bbb Q}_{\langle k \rangle}(p))) $.
In the notation of (3.1) we have
$$
\widetilde{N}^{i}H^{2i}(X,{\Bbb Q}_{\langle k \rangle}(i)) =
\Hdg^{i}(X,{\Bbb Q}_{\langle k \rangle})
\otimes {\Bbb Q}_{\langle k \rangle}.
$$

For a smooth connected curve
$ S $, let
$ {\Cal M}_{S} $ be the pull-back of
$ {\Cal M} $ by the structure morphism
$ a_{S} : S \rightarrow \Spec {\Bbb C} $.
Then we have a short exact sequence
$$
0 \rightarrow \Ext^{1}({\Bbb Q}_{\langle k \rangle}, {\Cal M})
\rightarrow \Ext^{1}({\Bbb Q}_{S,\langle k \rangle}, {\Cal M}_{S})
\rightarrow \Hom({\Bbb Q}_{\langle k \rangle},
H^{1}(a_{S})_{*}{\Cal M}_{S}) \rightarrow 0
$$
by the Leray spectral sequence for
$ (a_{S})_{*}{\Cal M}_{S} $, where the first morphism is induced by the
pull-back by
$ a_{S} $.
So the assertion is reduced to the vanishing of the last term (using a
stratification of the Hilbert scheme).
We have a natural isomorphism
$$
\Hom({\Bbb Q}_{\langle k \rangle},H^{1}(a_{S})_{*}{\Cal M}_{S}) =
\Hom({\Bbb D}{\Cal M}, H^{1}(S,{\Bbb Q}_{\langle k \rangle})),
$$
because
$ H^{1}(a_{S})_{*}{\Cal M}_{S} = {\Cal M}\otimes H^{1}(S,{\Bbb Q}_{\langle
k \rangle}) $.
(Here
$ {\Bbb D}{\Cal M} $ denotes the dual of
$ {\Cal M} $.)
Since
$ {\Bbb D}{\Cal M} $ is pure of weight
$ 2 $, and
$ {\Gr}_{2}^{W}H^{1}(S,{\Bbb Q}_{\langle k \rangle}) $ is isomorphic to a
direct sum of copies of
$ {\Bbb Q}_{\langle k \rangle}(-1) $, the assertion is further reduced to
$$
\Hom({\Bbb Q}_{\langle k \rangle}(1), {\Cal M}) = 0.
$$
But this is clear by the definition of
$ {\Cal M} $ (together with the semisimplicity of pure objects).
\enddemo

\medskip\noindent
{\bf 5.10.}
{\it Remark.} It is conjectured by C.~Voisin [47] that
$ {\CH}_{\ind}^{2}(X,1)_{\Bbb Q} $ is countable.
By (5.3), (5.7) and (5.9), this conjecture can be reduced to the
injectivity of (5.1.2) (or more precisely, to the hypothesis of (5.3)).

\proclaim{{\bf 5.11.~Proposition}}
Let
$ X $ be a smooth complex projective variety.
Assume
$ k $ is a number field.
Then the morphism (5.7.1) is surjective if
$ H^{2}(X,{\Bbb Q}_{\langle k \rangle})
/N^{1}H^{2}(X,{\Bbb Q}_{\langle k \rangle}) $ is global section-free in the
sense of (2.4) and if the Abel-Jacobi map in (0.4) is injective for
codimension
$ 2 $ cycles on any
$ k $-smooth projective models of
$ X $ (where the notion of
$ k $-smooth projective model is defined as in (0.4)).
\endproclaim

\demo\nofrills {Proof.\usualspace}
Let
$ \pi : X_{R} \rightarrow S = \Spec R $ be as in (2.1), and
$ \pi ' : X_{R'} \rightarrow S' = \Spec R' $ its base change by a finitely
generated smooth
$ k $-subalgebra
$ R' $ of
$ {\Bbb C} $ containing
$ R $.
By definition, the target of (5.7.1) is the inductive limit of
$$
\Ext^{1}({\Bbb Q}_{S',k_{R'}}, R^{2}\pi '_{*}{\Bbb Q}_{X_{R'}/k_{R'}}(2)
/N^{1}R^{2}\pi '_{*}{\Bbb Q}_{X_{R'}/k_{R'}}(2)),
$$
where
$ N^{1}R^{2}\pi '_{*}{\Bbb Q}_{X_{R'}/k_{R'}}(2) $ is isomorphic to a
direct sum of copies of
$ {\Bbb Q}_{S',k_{R'}}(1) $ if
$ R $ is sufficiently large.
This extension group is isomorphic to
$$
\Hom({\Bbb Q}_{k_{R'}}, H^{1}(S'/k_{R'}, R^{2}\pi '_{*}{\Bbb
Q}_{X_{R'}/k_{R'}}(2)
/N^{1}R^{2}\pi '_{*}{\Bbb Q}_{X_{R'}/k_{R'}}(2)))
$$
due to the first assumption (by using the Leray spectral sequence).
So the assertion is reduced to the following:
\enddemo

\proclaim{{\bf 5.12.~Proposition}}
Let
$ X $ be a smooth projective variety over a number field
$ k $ such that
$ k $ is algebraically closed in the function field of
$ X $.
Let
$ U $ be a dense open subvariety of
$ X $ such that
$ Y := X \backslash U $ is a divisor with normal crossings whose
irreducible components are smooth.
Then the morphism induced by the cycle map
$$
\CH^{2}(U,1)_{\Bbb Q} \rightarrow \Hom({\Bbb Q}_{k}, H^{3}(U /k,{\Bbb
Q}(2)))
$$
is surjective, if the Abel-Jacobi map in (0.4) is injective for codimension
$ 2 $ cycles on
$ X $.
\endproclaim

\demo\nofrills {Proof.\usualspace}
Let
$ M = H^{3}(U/k,{\Bbb Q}(2)) $ as in the proof of (4.7), and
$ \CH^{1}(Y{)}_{\Bbb Q}^{(0)} $ the kernel of the composition
$ \CH^{1}(Y)_{\Bbb Q} \rightarrow \CH^{2}(X)_{\Bbb Q} \rightarrow
H^{4}(X/k, {\Bbb Q}(2)) $.
Then we have a commutative diagram induced by the cycle map
$$
\CD
\CH^{2}(U,1)_{\Bbb Q} @>>> \CH^{1}(Y{)}_{\Bbb Q}^{(0)} @>>>
L^{1}\CH^{1}(X)_{\Bbb Q}
\\
@VVV @VVV @VVV
\\
\Hom({\Bbb Q}_{k}, M) @>>> \Hom({\Bbb Q}_{k}, M/W_{-1}M)
@>>> \Ext^{1}({\Bbb Q}_{k}, {\Gr}_{-1}^{W}M)
\endCD
$$
For
$ \xi \in \Hom({\Bbb Q}_{k}, M) $, let
$ \xi ' $ be its image in
$ \Hom({\Bbb Q}_{k}, M/W_{-1}M) $.
Then by the arguments in the proof of (4.7), there exists
$ \zeta ' \in \CH^{1}(Y{)}_{\Bbb Q}^{(0)} $ such that
$ \xi ' $ coincides with the image of
$ \zeta ' $ in
$ \Hom({\Bbb Q}_{k}, M/W_{-1}M) $.

Since
$ \Hom({\Bbb Q}_{k}, W_{-1}M) = 0 $, the first morphism in the bottom row
is injective.
So it is enough to show that we can choose
$ \zeta ' $ coming from
$ \CH^{2}(U,1)_{\Bbb Q} $.
By assumption, it is equivalent to that the image of
$ \zeta ' $ in
$ J^{2}(X/k)_{\Bbb Q} $ vanishes in the notation of (1.3.4)
(replacing
$ \zeta ' $ if necessary).
But the image of
$ \zeta ' $ in
$ \Ext^{1}({\Bbb Q}_{k}, {\Gr}_{-1}^{W}M) $ is zero, and hence
$ cl'(\zeta ') \in J^{2}(X/k)_{\Bbb Q} $ belongs to the image of
$ J^{1}(\widetilde{Y}/k)_{\Bbb Q} $, where
$ \widetilde{Y} $ denotes the normalization of
$ Y $.
So the assertion follows from (4.8).
This completes the proofs of (5.12) and (5.11).
\enddemo

\bigskip\bigskip
\centerline{{\bf References}}

\bigskip

\item{[1]}
Asakura, M., Motives and algebraic de Rham cohomology, preprint.

\item{[2]}
Beilinson, A., Higher regulators and values of
$ L $-functions, J. Soviet Math. 30 (1985), 2036--2070.

\item{[3]}
\SameAuthor, Notes on absolute Hodge cohomology, Contemporary Math. 55
(1986) 35--68.

\item{[4]}
\SameAuthor, Height pairing between algebraic cycles, Lect. Notes in Math.,
vol. 1289, Springer, Berlin, 1987, pp. 1--26.

\item{[5]}
\SameAuthor, On the derived category of perverse sheaves, ibid. pp. 27--41.

\item{[6]}
Beilinson, A., Bernstein, J. and Deligne, P., Faisceaux pervers,
Ast\'erisque, vol. 100, Soc. Math. France, Paris, 1982.

\item{[7]}
Bloch, S., Lectures on algebraic cycles, Duke University Mathematical
series 4, Durham, 1980.

\item{[8]}
\SameAuthor, Algebraic cycles and higher $ K $-theory, Advances in Math.,
61 (1986), 267--304.

\item{[9]}
\SameAuthor, Algebraic cycles and the Beilinson conjectures, Contemporary
Math. 58 (1) (1986), 65--79.

\item{[10]}
Bloch, S. and Srinivas, V., Remarks on correspondences and algebraic
cycles, Amer. J. Math. 105 (1983), 1235--1253.

\item{[11]}
Carlson, J., Extensions of mixed Hodge structures, in Journ\'ees de
G\'eom\'etrie Alg\'ebrique d'Angers 1979, Sijthoff -Noordhoff Alphen a/d
Rijn, 1980, pp. 107--128.

\item{[12]}
Deligne, P., Th\'eorie de Hodge I, Actes Congr\`es Intern. Math., 1970,
vol. 1, 425-430; II, Publ. Math. IHES, 40 (1971), 5--57; III ibid., 44
(1974), 5--77.

\item{[13]}
\SameAuthor, Th\'eor\`eme de finitude en cohomologie {\it l}-adique, in SGA
4 1/2, Lect. Notes in Math., vol. 569, Springer, Berlin, 1977, 233--261.

\item{[14]}
\SameAuthor, Valeurs de fonctions $ L $ et p\'eriodes d'int\'egrales, in
Proc. Symp. in pure Math., 33 (1979) part 2, 313--346.

\item{[15]}
Deligne, P., Milne, J., Ogus, A. and Shih, K., Hodge Cycles, Motives, and
Shimura varieties, Lect. Notes in Math., vol 900, Springer, Berlin, 1982.

\item{[16]}
Deninger, C. and Scholl, A., The Beilinson conjectures, in Proceedings
Cambridge Math. Soc. (eds. Coats and Taylor) 153 (1992), 173--209.

\item{[17]}
El Zein, F. and Zucker, S., Extendability of normal functions associated to
algebraic cycles, in Topics in transcendental algebraic geometry, Ann.
Math. Stud., 106, Princeton Univ. Press, Princeton, N.J., 1984, pp.
269--288.

\item{[18]}
Esnault, H. and Levine, M., Surjectivity of cycle maps, Ast\'erisque 218
(1993), 203--226.

\item{[19]}
Fulton, W., Intersection theory, Springer, Berlin, 1984.

\item{[20]}
Gordon, B. and Lewis, J., Indecomposable higher Chow cycles on products of
elliptic curves, J. Alg. Geom. 8 (1999), 543--567.

\item{[21]}
Green, M., Griffiths' infinitesimal invariant and the Abel-Jacobi map, J.
Diff. Geom. 29 (1989), 545--555.

\item{[22]}
\SameAuthor, What comes after the Abel-Jacobi map?, preprint.

\item{[23]}
\SameAuthor, Algebraic cycles and Hodge theory (Lecture notes at the Banff
conference, based on a collaboration with P. Griffiths).

\item{[24]}
Griffiths, P., On the period of certain rational integrals I, II, Ann.
Math. 90 (1969), 460--541.

\item{[25]}
Jannsen, U., Mixed motives and algebraic $ K $-theory, Lect. Notes in
Math., vol. 1400, Springer, Berlin, 1990.

\item{[26]}
\SameAuthor, Motivic sheaves and filtrations on Chow groups, Proc. Symp.
Pure Math. 55 (1994), Part 1, pp. 245--302.

\item{[27]}
Kashiwara, M., A study of variation of mixed Hodge structure, Publ. RIMS,
Kyoto Univ. 22 (1986) 991--1024.

\item{[28]}
Kleiman, S., Algebraic cycles and Weil conjecture, in Dix expos\'es sur la
cohomologie des sch\'emas, North-Holland, Amsterdam, 1968, pp. 359 --386.

\item{[29]}
Lang, S., Abelian varieties, Wiley, New York, 1959.

\item{[30]}
M\"uller-Stach, S., Constructing indecomposable motivic cohomology classes
on algebraic surfaces, J. Alg. Geom. 6 (1997), 513--543.

\item{[31]}
Mumford, D., Rational equivalence of
$ 0 $-cycles on surfaces, J. Math. Kyoto Univ. 9 (1969), 195--204.

\item{[32]}
Murre, J., On the motive of an algebraic surface, J. Reine Angew. Math. 409
(1990), 190--204.

\item{[33]}
\SameAuthor, On a conjectural filtration on Chow groups of an algebraic
variety, Indag. Math. 4 (1993), 177--201.

\item{[34]}
Roitman, A., Rational equivalence of zero cycles, Math. USSR Sbornik 18
(1972), 571--588.

\item{[35]}
Saito, M., Modules de Hodge polarisables, Publ. RIMS, Kyoto Univ. 24 (1988)
849--995.

\item{[36]}
\SameAuthor, Mixed Hodge Modules, Publ. RIMS, Kyoto Univ., 26 (1990),
221--333.

\item{[37]}
\SameAuthor, On the formalism of mixed sheaves, RIMS-preprint 784, Aug.
1991.

\item{[38]}
\SameAuthor, Admissible normal functions, J. Alg. Geom. (1996), 235--276.

\item{[39]}
\SameAuthor, Tate conjecture and mixed perverse sheaves, preprint.

\item{[40]}
\SameAuthor, Hodge conjecture and mixed motives, I, Proc. Symp. Pure Math.
53 (1991), 283--303; II, in Lect. Notes in Math., vol. 1479, Springer,
Berlin, 1991, pp. 196--215.

\item{[41]}
\SameAuthor, Some remarks on the Hodge type conjecture, Proc. Symp. Pure
Math. 55 (1991), Part 1, pp. 85--100.

\item{[42]}
\SameAuthor, Bloch's conjecture, Deligne cohomology and higher Chow group,
preprint.

\item{[43]}
Saito, S., Motives and filtrations on Chow groups, Inv. Math. 125 (1996),
149 --196.

\item{[44]}
Steenbrink, J. and Zucker, S., Variation of mixed Hodge structure I, Inv.
Math. 80 (1985) 489--542.

\item{[45]}
Voisin, C., Variations de structures de Hodge et z\'ero-cycles sur les
surfaces g\'en\'erales, Math. Ann. 299 (1994), 77--103.

\item{[46]}
\SameAuthor, Transcendental methods in the study of algebraic cycles, in
Lect. Notes in Math. vol. 1594, pp. 153--222.

\item{[47]}
\SameAuthor, Remarks on zero-cycles of self-products of varieties, in
Moduli of Vector Bundles, Lect. Notes in Pure and Applied Mathematics, vol.
179, M. Dekker, New York, 1996, pp. 265--285.

\item{[48]}
\SameAuthor, Some results on Green's higher Abel-Jacobi map, Ann. Math. 149
(1999), 451--473.

\bigskip
\noindent
Nov. 24, 2000
\bye